\documentclass{gtart}

\def\ifplaintex{\expandafter\ifx\csname documentclass\endcsname\relax}

\def\gtp{{\mathsurround=0pt\it $\cal G\mskip-2mu$eometry \&\ 
$\cal T\!\!$opology $\cal P\!$ublications}}  

\def\recd{{\small Received:\qua\receiveddate\ifx\reviseddate\relax
\else\qquad Revised:\qua\reviseddate\fi\par}} 

\def\lognumber#1{\def\thelognumber{#1}}
\def\volumenumber#1{\def\thevolumenumber{#1}}
\def\volumeyear#1{\def\thevolumeyear{#1}}
\def\papernumber#1{\def\thepapernumber{#1}}
\def\pagenumbers#1#2{\def\startpage{#1}\def\finishpage{#2}}
\def\published#1{\def\publishdate{#1}}

\def\received#1{\def\receiveddate{#1}}
\def\revised#1{\def\reviseddate{#1}}
\def\accepted#1{\def\accepteddate{#1}}

\def\asciiemail#1{\def\theasciiemail{#1}}

\long\def\asciiabstract#1{\long\def\theasciiabstract{#1}}

\let\\\par\let\thelognumber\relax\let\thevolumenumber\relax
\let\thepapernumber\relax\let\thevolumeyear\relax\let\startpage\relax
\let\finishpage\relax\let\publishdate\relax\let\receiveddate\relax
\let\reviseddate\relax\let\accepteddate\relax\let\theasciititle\relax
\let\theasciiauthors\relax
\let\theasciiabstract\relax

\let\theasciiemail\relax

\ifplaintex
\font\logobig=cmssbx10 scaled 3836
\font\logomed=cmssbx10 scaled 2557
\else
\font\logobig=cmssbx10 scaled 4200
\font\logomed=cmssbx10 scaled 2800
\fi

\long\def\makeagttitle{   
\count0=\startpage
\agt\hfill      
\hbox to 45truept{\vbox to 0pt{\vglue -13truept{\logomed A\kern -.37em{\logobig 
T}\kern -.38em G}\vss}\hss}
\break
{\small Volume \thevolumenumber\ (\thevolumeyear)
\startpage--\finishpage\nl
Published: \publishdate}

\vglue .25truein

{\parskip=0pt\leftskip 0pt plus
1fil\def\\{\par\smallskip}{\Large\bf\thetitle}\par\medskip} \vglue
0.05truein

{\parskip=0pt\leftskip 0pt plus 1fil\def\\{\par}{\sc\theauthors}
\par\medskip}%
 
\vglue 0.03truein 

{\small\leftskip 25truept\rightskip 25truept{\bf Abstract}\stdspace\theabstract

{\bf AMS Classification}\stdspace\theprimaryclass
\ifx\thesecondaryclass\relax\else; \thesecondaryclass\fi\par
{\bf Keywords}\stdspace \thekeywords\par}\vglue 7truept

}   

\ifplaintex
\font\phead=cmsl9 scaled 950
\font\pnum=cmbx10 scaled 913
\font\pfoot=cmsl9 scaled 950
\headline{\vbox to 0pt{\vskip -4.5mm\line{\small\phead\ifnum
\count0=\startpage ISSN 1472-2739 (on-line) 1472-2747 (printed)
\hfill {\pnum\folio}\else\ifodd\count0\def\\{ }%
\ifx\theshorttitle\relax\thetitle\else\theshorttitle\fi\hfill{\pnum\folio}
\else\def\\{ and }{\pnum\folio}\hfill\ifx\theshortauthors\relax\theauthors
\else\theshortauthors\fi\fi\fi}\vss}}
\footline{\vbox to 0pt{\vglue 0mm\line{\small\pfoot\ifnum\count0=\startpage
\copyright\ \gtp\hfill\else
\agt, Volume \thevolumenumber\ (\thevolumeyear)\hfill\fi}\vss}}
\else
\font=cmsl9 scaled 1050
\font\lnum=cmbx10 
\font=cmsl9 scaled 1050
\makeatletter
\def\@oddhead{{\small\lhead\ifnum\count0=\startpage ISSN 1472-2739 
(on-line) 1472-2747 (printed)\hfill {\lnum\number\count0}\else\ifodd\count0
\def\\{ }\ifx\theshorttitle\relax \thetitle \else\theshorttitle\fi\hfill
{\lnum\number\count0}\else\def\\{ and }{\lnum\number\count0}
\hfill\ifx\theshortauthors\relax 
\theauthors\else\theshortauthors\fi\fi\fi}}\def\@evenhead{\@oddhead}
\def\@oddfoot{\small\lfoot\ifnum\count0=\startpage\copyright\ \gtp\hfill\else
\agt, Volume \thevolumenumber\ (\thevolumeyear)\hfill\fi}
\def\@evenfoot{\@oddfoot}
\makeatother
\fi
\let\maketitlepage\makeagttitle
\let\makeshorttitle\maketitlepage
\let\maketitle\maketitlepage

\newwrite\gtoutfile
\long\gdef\makeheadfile{  
{\def\\{, }\def\s{ }
\immediate\openout\gtoutfile head.xxx
\immediate\write\gtoutfile{To: math@arxiv.org}
\immediate\write\gtoutfile{Subject: put OR rep NNNNN:ppppp}
\immediate\write\gtoutfile{--text follows this line--}
\immediate\write\gtoutfile{Proxy-for: \ifx\theasciiauthors\relax
\theauthors\else\theasciiauthors\fi\s<\ifx\theasciiemail\relax\theemail\else\theasciiemail\fi>}
\immediate\write\gtoutfile{\noexpand\\}
\immediate\write\gtoutfile{Authors: \ifx\theasciiauthors\relax
\theauthors\else\theasciiauthors\fi}
{\def\\{ }\immediate\write\gtoutfile{Title: \ifx\theasciititle\relax
\thetitle\else\theasciititle\fi}}
\immediate\write\gtoutfile{Subj-class: GT or SG, GR etc}
\immediate\write\gtoutfile{MSC-class: \theprimaryclass\ifx\thesecondaryclass\relax\else, \thesecondaryclass\fi}
\immediate\write\gtoutfile{Journal-ref: Algebr. Geom. Topol. \thevolumenumber\s
(\thevolumeyear) \startpage-\finishpage}
\immediate\write\gtoutfile{Comments: Published by Algebraic and
Geometric Topology at}
\immediate\write\gtoutfile{\s\s\s  http://www.maths.warwick.ac.uk/agt/AGTVol\thevolumenumber/agt-\thevolumenumber-\thepapernumber.abs.html}
\immediate\write\gtoutfile{\noexpand\\}
\immediate\write\gtoutfile{}
\ifx\theasciiabstract\relax
\immediate\write\gtoutfile{\theabstract}\else
\immediate\write\gtoutfile{\theasciiabstract}\fi
\immediate\write\gtoutfile{}
\immediate\write\gtoutfile{\noexpand\\}
\immediate\write\gtoutfile{}
\immediate\closeout\gtoutfile}}  

\def\maketitlepage{\makeagttitle\makeheadfile}
\let\makeshorttitle\maketitlepage
\let\maketitle\maketitlepage

\lognumber{14}
\volumenumber{3}
\volumeyear{2003}
\papernumber{14}
\published{22 May 2003}
\pagenumbers{435}{472}
\received{9 August 2002}
\revised{7 February 2003}
\accepted{4 April 2003}

\usepackage{amsmath,graphicx,amssymb}

\newtheorem{pro}{Proposition}[section]
\newtheorem{thm}[pro]{Theorem}
\newtheorem{lem}[pro]{Lemma}
\newtheorem{clm}{Claim}

\newtheorem{cor}[pro]{Corollary}
\newtheorem*{begthm}{Theorem}
\theoremstyle{definition}
\newtheorem{dfn}{Definition}[section]
\def\strut{\vrule width0pt height 12pt depth7pt}

\theoremstyle{remark}
\newtheorem*{rems}{Remarks}

\newtheorem*{note}{Note}

\begin{document}

\title{Small Seifert-fibered Dehn surgery\\on hyperbolic knots}

\author{John C. Dean}
\address{818 Forest Ave., Oak Park, IL 60302, USA} 
\email{dean\_john@yahoo.com} 
\asciiemail{dean_john@yahoo.com}


\begin{abstract} In this paper, we define the primitive/Seifert-fibered property for a knot in $S^3$.  If satisfied, the property ensures that the knot has a Dehn surgery that yields a small Seifert-fibered space (i.e.\ base $S^2$ and three or fewer critical fibers).  Next we describe the twisted torus knots, which provide an abundance of examples of primitive/Seifert-fibered knots.  By analyzing the twisted torus knots, we prove that nearly all possible triples of multiplicities of the critical fibers arise via Dehn surgery on primitive/Seifert-fibered knots.
\end{abstract}

\asciiabstract{
In this paper, we define the primitive/Seifert-fibered property for a
knot in S^3.  If satisfied, the property ensures that the knot has a
Dehn surgery that yields a small Seifert-fibered space (i.e. base
S^2 and three or fewer critical fibers).  Next we describe the
twisted torus knots, which provide an abundance of examples of
primitive/Seifert-fibered knots.  By analyzing the twisted torus
knots, we prove that nearly all possible triples of multiplicities of
the critical fibers arise via Dehn surgery on
primitive/Seifert-fibered knots.}

\maketitle

\section{Introduction}
 
Since all compact orientable 3-manifolds can be realized by Dehn
surgery on a link in the 3-sphere \cite{lickorish:1962, wallace:1960},
considerable energy has been devoted to trying to understand Dehn
surgery.  Thurston showed that hyperbolic knots are ubiquitous.
Moreover, he showed that only finitely many Dehn surgeries on a
hyperbolic knot can be non-hyperbolic.  For this reason, a
non-hyperbolic Dehn surgery on a hyperbolic knot is called an
\emph{exceptional surgery}.

A \emph{small Seifert-fibered space} is a Seifert-fibered space with base space $S^2$ and at most three singular fibers.  If Thurston's geometrization conjecture for 3-manifolds is true, then any exceptional surgery must either contain an essential sphere or torus, or be a small Seifert-fibered space.  Much is known about toroidal exceptional surgeries, while the cabling conjecture \cite{g-as:1986} would imply that no hyperbolic knot has a Dehn surgery manifold with an essential sphere.

Until recently, relatively little was known about small Seifert-fibered exceptional surgeries except when the surgered manifold has finite fundamental group.  In this case, known results include the knot complement theorem and the cyclic surgery theorem \cite{gl:1989, cgls:1987}.  Boyer and Zhang have proved theorems which limit the behavior of the slopes of more general types of finite surgeries \cite{bz:1996}.

One result on exceptional surgeries that \emph{does} apply to small Seifert-fibered Dehn surgeries with infinite fundamental group is the $2\pi$ theorem of Gromov and Thurston (see \cite{bh:1996}).  The theorem implies that there can be at most 24 non-negatively curved Dehn surgeries on a hyperbolic knot.  Since Seifert-fibered spaces are known not to be negatively curved, there can be at most 24 Seifert-fibered Dehn surgeries on a hyperbolic knot.  Recently Agol \cite{agol:1998} and Lackenby \cite{lackenby:2000} have improved the $2\pi$ theorem.  In particular, Agol showed that there are at most twelve surgeries on a hyperbolic knot that are either reducible, toroidal, or Seifert-fibered.

The first known examples of small Seifert manifolds arising from Dehn surgery on hyperbolic knots were given by \cite{fs:1980}.  Berge has a construction which produces families of knots with lens space Dehn surgeries \cite{berge:u1}. Many of these knots are hyperbolic.  It is an open question whether or not the Berge knots include all knots with lens space Dehn surgeries.  Berge has explicitly described the knots which arise from his construction.

Various examples of small Seifert-fibered Dehn surgeries have been given by Bleiler and Hodgson \cite{bh:1996}, Eudave-Mu\~noz \cite{munoz:1993}, Boyer and Zhang \cite{bz:1996}, and Miyazaki and Motegi \cite{mm:1999}.

The Bleiler-Hodgson and Boyer-Zhang examples arose in trying to understand exceptional surgeries with finite fundamental groups.  The examples of Eudave-Mu\~noz were constructed as examples of hyperbolic knots with toroidal $\mathbb{Z} /2$ surgeries.  Eudave-Mu\~noz showed that these knots often have small Seifert-fibered surgeries as well.

The work in this paper is the author's thesis, together with generalizations and improvements of those results.  Since then much interesting and significant related work has been completed, in particular \cite{munoz:u2} and \cite {mmm:u1}. 

The author would like to thank his thesis adviser,  Cameron Gordon, for many helpful suggestions and enlightening conversations as this work progressed.

\subsection{Hyperbolic Knots with small Seifert-fibered surgeries}

In this paper we describe a new construction of knots with small Seifert-fibered Dehn surgeries.  Knots that arise from the construction are called primitive/Seifert-fibered knots. The construction is a generalization of Berge's construction of knots with lens space Dehn surgeries.  We will sometimes abbreviate primitive/Seifert-fibered as P/SF; and use SSFS for small Seifert-fibered space.

In Section \ref{s:psf}, we describe the primitive/Seifert-fibered construction and we show that a primitive/Seifert-fibered knots is guaranteed to have a small Seifert-fibered Dehn surgery, or one that is the connected sum of two lens spaces.  According to the cabling conjecture, the latter possibility would not arise from any hyperbolic knot.  To show that non-trivial examples exist, we demonstrate how the slope $2/1$ and slope $3/1$ surgeries on the twist knots conform to the P/SF construction.

While there are other primitive/Seifert-fibered knots, we focus in this paper on those that are twisted torus knots, a notion that we define in Section \ref{s:ttor}.  While only certain twisted torus knots are primitive/Seifert-fibered, those that are provide a rich set of examples of the phenomenon, and the author has explored their structure in some detail.  Group theoretic properties of certain canonical embeddings of twisted torus knots on a genus two Heegaard surface are also developed in Section~\ref{s:ttor}.  As a corollary of this work, we show that any possible one-relator presentation of the group $\langle x, y \ | \ x^my^n \rangle$ can be realized geometrically (i.e.\ by adding a 2-handle to a handlebody of genus two).

These properties are exploited in Section \ref{s:class&mult} to classify those P/SF twisted torus knots that are middle Seifert-fibered, a notion that is defined in Section \ref{ss:sfwords}.  After classifying the primitive/middle-Seifert-fibered twisted torus knots, the multiplicities of the critical fibers are calculated for the SSFS arising from Dehn surgery for each such knot.

Finally, in Section \ref{s:ubiquity} we show that, as measured by the multiplicities of their critical fibers, many small Seifert manifolds can be realized by Dehn surgery on these knots. In particular, we have the following result.
\begin{begthm}
For any triple of integers $(\mu_1,\mu_2,\mu_3)$ with $\gcd(\mu_1, \mu_2)= 1$, there is a non-torus knot with a small Seifert-fibered Dehn surgery with these multiplicities.
\end{begthm}
Many of the knots (described in section 5) used to prove this theorem  are known to be hyperbolic, and we expect that all of them are hyperbolic.  A similar result appeared in \cite{dean:1996} with the added restriction that $|\mu_1-~\mu_2| > 1$.  The improved result has also been obtained by Miyazaki and Motegi using the examples in \cite{mm:1999} discussed below.  Moreover, Miyazaki and Motegi have shown their knots to be hyperbolic.

Of the known example of knots with SSF surgeries, many have been shown to satisfy the P/SF construction (including those in (\cite{munoz:1993} and \cite{mm:1999}).  Recently, Mattman, Mizaki, and Motegi showed \cite{mmm:u1} that there is a hyperbolic knot with a small Seifert-fibered Dehn surgery that does not arise via the primitive/Seifert-fibered construction.

\section{Primitive/Seifert-fibered Knots}

\label{s:psf}
In this section we will describe a way to construct knots in $S^3$ that have a Dehn surgery that is a Seifert-fibered space with base $S^2$ and three or fewer critical fibers.  The construction is a generalization of Berge's construction of knots with lens space Dehn surgeries \cite{berge:u1}.  From the definition of the construction, it is not clear that any nontrivial non-Berge examples arise, hence we will describe a simple family of nontrivial examples that arise from the construction.  We begin with the definitions of some relevant concepts.  We will consider only orientable 3-manifolds throughout.

\subsection{Knots in separating surfaces,\,the surface slope,\,and\,2-handle addition}

We begin by defining the notion of 2-handle addition for a 3-manifold with boundary.

\begin{dfn}
\label{d:2handle}
Let $\gamma$ be a simple closed curve in the boundary of a 3-manifold $M$, and let $A$ be a regular neighborhood of $\gamma$ in $\partial M$.  Then $M \cup_A (D^2 \times I)$ is the result of \emph{2-handle addition along $\gamma$}, where $A$ and $\partial D^2 \times I$ are identified.
\end{dfn}

Next we define the surface slope for a knot contained in a surface in a 3-manifold, and show how Dehn surgery at this slope is related to 2-handle addition when the surface is separating.

\noindent
{\bf Notation}\qua  $M(K, \gamma)$ denotes the manifold obtained from Dehn surgery on $K$ with slope $\gamma$.  $N(K)$ is a regular neighborhood of a knot $K$.

\begin{dfn}
\label{d:surfslope}
If K is a knot embedded in a surface $F$ in a 3-manifold then the isotopy class in $\partial N(K)$ of the curve(s) in $\partial N(K) \cap F$ is called the \emph{surface slope of $K$ with respect to $F$.} 
\end{dfn}

\begin{lem}
\label{l:surfslope}
Let K be a knot contained in a separating surface $F$ in a 3-manifold $M$, i.e.  $M=V\cup_FV'$, and let $\gamma$ be the surface slope of $K$ with respect to $F$.  Then $M(K, \gamma) \cong W\cup_{\tilde F}W'$ where W (resp.\  W') is obtained from V (resp.\ V') by attaching a 2-handle along $K$, and $\tilde F=(F-N(K)) \cup (D^2 \times \{ 0, 1 \})$. 
\end{lem}

\proof
Let $A$ (resp.\ $A'$) be the annulus $\partial(N(K)) \cap V$ (resp.\ $\partial(N(K)) \cap V'$),and let $c_1$ and $c_2$ be their (shared) boundary curves. 
Denote the Dehn surgery solid torus by $U$, and let $h : \partial U \rightarrow \partial(S^3 - N(K))$ be the attaching map for the Dehn surgery.
\begin{figure}[ht!]
\begin{center}
\includegraphics[bb=60 258 468 641,width=5in]{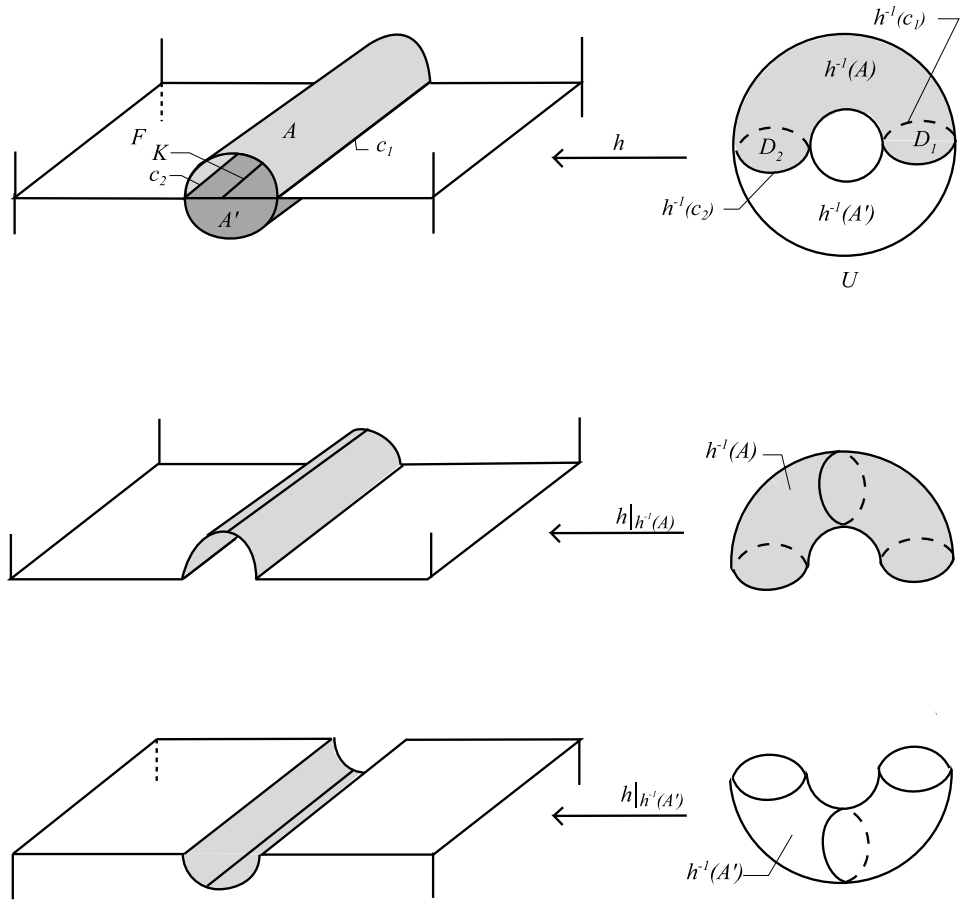}
\caption{Surface slope Dehn surgery}
\label{f:surf_slope}
\end{center}
\end{figure}

Since we are considering surface slope Dehn surgery, the curves $h^{-1}(c_1)$ and $h^{-1}(c_2)$ bound disks $D_1$ and $D_2$ in $U$.  We may cut the Dehn surgered manifold along $(F-N(K)) \cup D_1 \cup D_2$.  The resulting pieces (see Figure \ref{f:surf_slope}) are homeomorphic to $$W = V \cup_{K} \text{ 2-handle}$$
and $$W' = V' \cup_{K} \text{ 2-handle}.\eqno{\qed}$$

\subsection{Special elements in a free group and curves on handlebodies}

Let $G_{a,b}$ denote the group $\langle x, y | x^ay^b \rangle$.  When $a$ and $b$ are coprime this group is the fundamental group of the $(a,~b)$ torus knot.  More generally, $G_{a, b}$ is the fundamental group of a Seifert-fibered space over the disk with two critical fibers of multiplicity $a$ and $b$.  Recall that a \emph{basis} for a free group is a set of elements that freely generates the group.

\begin{dfn}
\label{d:prim}
An element in a free group is \emph{primitive} if it is part of a basis.
\end{dfn}
\begin{dfn}
\label{d:sf}
An element $w$ in the free group on $x$ and $y$ is \emph{($a$,~$b$) Seifert-fibered} if $\langle x, y | w \rangle \cong G_{a,b}$ for integers $a$ and $b$ both non-zero.  
\end{dfn}
Note that a primitive element in $\langle x, y \rangle$ is Seifert-fibered.

Let $\gamma$ be a simple closed curve contained in the boundary of a genus two handlebody $H$. Since $\gamma$ represents an element (defined up to conjugacy) of $\pi_1(H)$, which is a free group of rank two, we say that $\gamma$ is \emph{primitive with respect to $H$} if it represents a primitive element in $\pi_1(H)$.  We define Seifert-fibered simple closed curves on the boundary of a genus two handlebody similarly.

The following lemma establishes the link between Seifert-fibered curves on a genus two handlebody and Seifert-fibered spaces.

\begin{lem} 
\label{l:sfhandle}
Let $\gamma$ be a curve in the boundary of a genus two handlebody $H$ that is $(a,b)$ Seifert-fibered with respect to $H$.  Then the manifold $M$ obtained by adding a 2-handle to $H$ along $\gamma$ is a Seifert-fibered space over $D^2$ with two critical fibers of multiplicities $a$ and $b$.  In particular, $$ M \cong D^2 \times S^1 \Leftrightarrow a \text{ or } b \text{ equals 1 } \Leftrightarrow \gamma \text{ is primitive.} $$
\end{lem}

\proof
The fundamental group of $M$ is  $G_{a,b}$, which has a non-trivial center. $M$ is irreducible and Haken, hence by \cite{waldhausen:1967}, is a Seifert-fibered manifold.  It is known that a Seifert-fibered manifold with such a fundamental group is a Seifert-fibered space with base space a disk and critical fibers of multiplicity $a$ and $b$.  By considering when $G_{a,b}$ is isomorphic to $\mathbb{Z}$, the last part follows.
\endproof

\subsection{Primitive/Seifert-fibered knots}

\label{ss:psf}
Putting the these definitions and lemmas together, we describe a property of a knot that ensures that the knot will have a Dehn surgery that is a small Seifert-fibered space or a connected sum of two lens spaces.

\begin{dfn} 
\label{d:psf}
Let $K$ be a knot contained in a genus two Heegaard surface $F$ for $S^3$, that is, $S^3 = H \cup_F H'$, where $H$ and $H'$ are genus two handlebodies.  Then $K$ is \emph{primitive/Seifert-fibered with respect to $F$} if it is primitive with respect to $H$ and Seifert-fibered with respect to $H'$.
\end{dfn}

\begin{pro}
\label{p:psf}
If a knot $K$ in $S^3$ is primitive/Seifert-fibered with respect to a genus two Heegaard surface, then Dehn surgery at the surface slope is either a small Seifert-fibered space or a connected sum of two lens spaces.  
\end{pro}

\proof
By Lemma \ref{l:surfslope} and Lemma \ref{l:sfhandle}, the Dehn surgered manifold is the union along a torus of a Seifert-fibered space over the disk with at most two critical fibers and a solid torus.  Thus surface slope Dehn surgery on a primitive/Seifert-fibered knot results in a manifold that is a Dehn filling of a Seifert-fibered space over the disk with two critical fibers.  

Since any non-meridinal simple closed curve on the boundary of a solid torus can be extended to a Seifert fibration of the solid torus, only one Dehn filling on such a Seifert-fibered manifold may fail to be Seifert-fibered.  This occurs exactly when the slope is an ordinary fiber.  For any other filling, a new critical fiber appears with multiplicity equal to the algebraic intersection number of the slope with the ordinary fiber.  So for any Dehn filling but one, the resulting manifold is a Seifert-fibered space over the sphere with at most three critical fibers, i.e.\ a small Seifert-fibered space.

A Seifert-fibered space over the disk with two fibers is the union of two solid tori glued along an annulus. When each fiber is non-trivial, this annulus intersects a meridian of each solid torus algebraically more than once.  A curve parallel to the annulus is an ordinary fiber.  When a solid torus is attached with slope equal to the ordinary fiber, the resulting manifold can be cut apart into two solid tori, each with a 2-handle attached along a curve which intersects the meridian more than once algebraically.  Thus each piece is a punctured lens space, so the manifold is a connected sum of two lens spaces.
\endproof

Boileau, Rost, and Zieschang have classified those curves $\gamma$ on the boundary of an abstract genus two handlebody $H$ that are Seifert-fibered \cite{brz:1988}.  An embedding of such a pair $(H, \gamma)$ into $S^3$ such that $H$ is unknotted and $\gamma$ is primitive with respect to $S^3 - H$ would give a P/SF knot.  However, it would be difficult to consider all possible unknotted embeddings of such pairs, and to determine which are primitive on the ``outside'' handlebody.

\rk{Remarks}
\begin{itemize}
\item According to the cabling conjecture, only cabled knots have reducible Dehn surgeries \cite{g-as:1986}.  In particular, the conjecture implies that a hyperbolic primitive/Seifert-fibered knot would always have a small Seifert-fibered surgery.

\item One could define a knot to be primitive/Seifert fibered in any 3-manifold of Heegaard genus less than or equal to two and the proposition would hold.

\item When the knot is primitive/primitive (\emph{doubly primitive}), a lens space results from the surface slope Dehn surgery (this is Berge's construction mentioned above).

\item Surface slope Dehn surgery on a doubly Seifert-fibered knot results in the union along a torus of two Seifert manifolds over the disk, each with two critical fibers.  Such a manifold is either a graph manifold or a Seifert manifold with base $S^2$ and four critical fibers.  No example is known of a hyperbolic knot with a Dehn surgery of the latter type. However, there are satellite knots with such Dehn surgeries \cite{mm:1997}.
\end{itemize}

\begin{figure}[ht!]
\begin{center}
\includegraphics[bb=98 480 471 554]{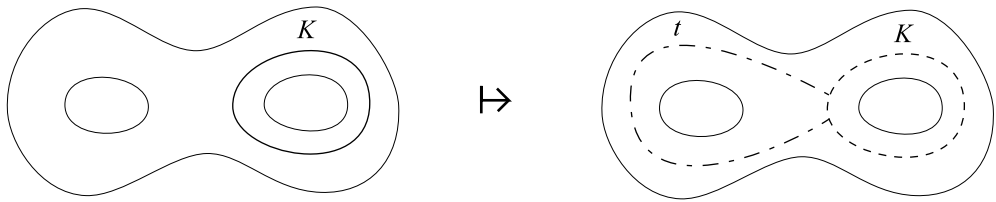}
\caption{Primitive on one side implies tunnel number one}
\label{f:tunnel}
\end{center}
\end{figure}

Note that any primitive/Seifert-fibered knot has tunnel number $1$.  In fact, any knot that is primitive with respect to one side of a genus two Heegaard surface has tunnel number $1$.  This is true since, by \cite{zieschang:1970}, there is a homeomorphism of the handlebody after which the knot $K$ appears as in Figure \ref{f:tunnel}. 
If one pushes $K$ slightly into the handlebody, removes a regular neighborhood of $K$, and then removes an appropriate tunnel $t$, then what remains in the handlebody is the product of a surface and an interval (see Figure \ref{f:tunnel}).  Thus the complement in $S^3$ is a handlebody, so the knot has tunnel number 1.

It is not clear from the definition of a primitive/Seifert-fibered knot whether any non-trivial examples exist.  We do not consider torus knots that arise from the construction to be interesting since Dehn surgery on torus knots is completely understood.  Berge's work shows that, in fact, there are a plethora of interesting knots that are doubly primitive, many of which are known to be hyperbolic.  In the next section we discuss a simple family of well-known hyperbolic knots, each of which has small Seifert-fibered Dehn surgeries at slopes $1$, $2$, and $3$.  We show that the slope $2$ and $3$ surgeries conform to the primitive/Seifert fibered construction.  

\subsection{The Twist Knots}

One of the simplest families of knots is the twist knots, which are obtained from the Whitehead link by performing a $1/n$ Dehn surgery on one of the components. In this family of knots $K_n$, indexed by the integers, only the unknot and the trefoil are non-hyperbolic.  The hyperbolic twist knots ($n \ne 0$,~$-1$) each have exceptional surgeries that are small Seifert-fibered spaces for the slopes 1, 2, and 3. Each twist knot has a toroidal surgery at slopes 0 and 4. The Figure eight is amphicheiral, hence it has small Seifert-fibered exceptional surgeries at all six slopes $\pm 1, \pm 2, \pm 3$, and toroidal exceptional surgeries at slopes 0 and $\pm 4$.  In \cite{bw:2001}, it is proved that these are the only exceptional surgeries for (hyperbolic) twist knots.  In fact, they show that the twist knots are the only two-bridge knots (except for the $(2, n)$ torus knots) that have any exceptional surgeries.

In this section we show that the slope $2$ and $3$ small Seifert-fibered surgeries on each twist knot can be realized as examples of the primitive/Seifert-fibered phenomenon.  We do not know whether or not the slope $1$ surgery can be realized by the primitive/Seifert-fibered construction.

\begin{figure}[ht!]
\begin{center}
\includegraphics[bb=110 309 500 675,width=5in]{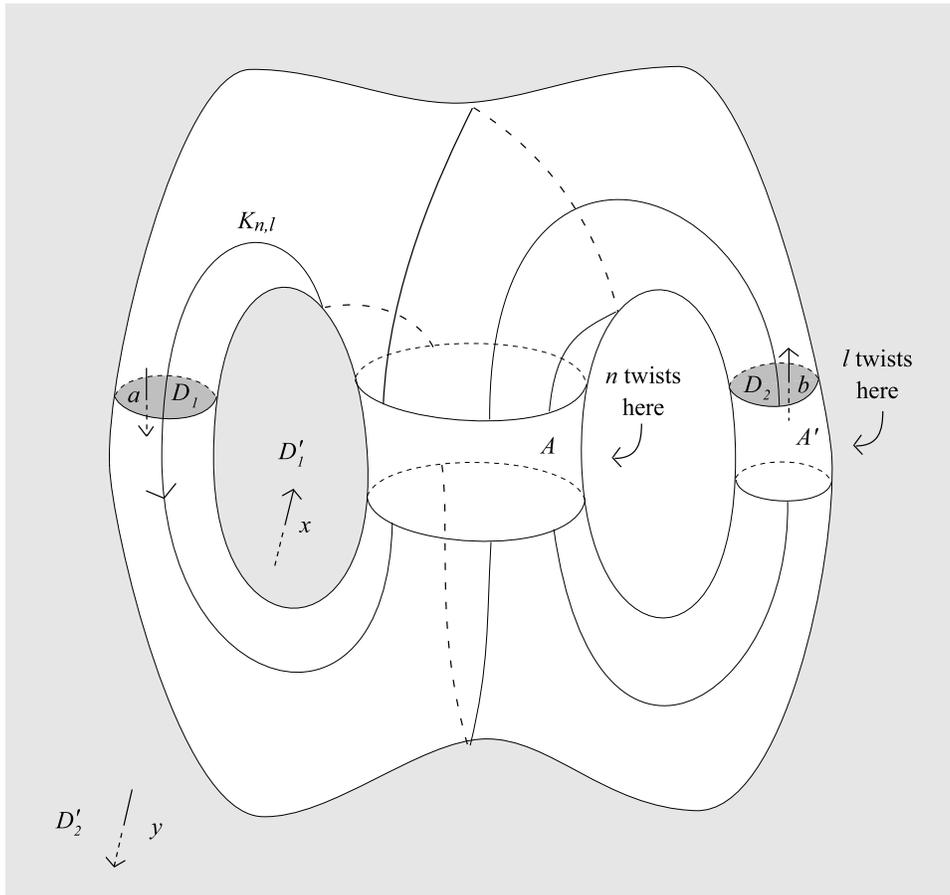}
\caption{Embeddings $K_{n,l}$ of the twist knots $K_n$}
\label{f:tw_knot}
\end{center}
\end{figure}

For slopes $2$ and $3$, and each integer $n$, we must exhibit an embedding of the twist knot $K_n$ on a genus two Heegaard surface $F$ in $S^3$ with this surface slope such that $K_n$ is primitive/Seifert-fibered with respect to $F$. The 2-bridge picture in Figure \ref{f:tw_knot} will allow us to do this.

Each twist knot $K_n$ can be realized by adding $n$ full twists to the annulus $A$.  Note that the way we have indexed the twist knots, $K_1$ is the figure-eight knot, while $K_{-1}$ is the right-handed trefoil.  Performing Dehn twists on the annulus $A'$ does not change the knot type of $K_n$, but it does change the surface slope.  Let $K_{n,l}$ denotes the embedding of $K_n$ in the Heegaard surface $F$ for $S^3$ with $l$ Dehn twists about $A'$.  A calculation shows that the surface slope of $K_{n,l}$ equals $2 + l$.  To calculate this surface slope, check that a full twist on the annulus $A$ contributes zero to the surface slope.

We will use the oriented disk system $D_1$, $D_2$ corresponding to generators $a$ and $b$ for $\pi_1(H)$, and $D_1', D_2'$ corresponding to generators $x$ and $y$ for $\pi_1(H')$.   Note that no matter what the value of $n$ or $l$, $K_{n, l}$ represents $ab$ in $\pi_1(H)$; hence it is primitive with respect to $H$.

Let $w_{n,l}$ be the conjugacy class of $K_{n,l}$ in $\pi_1(H')$.  Then a direct calculation shows that $w_{n,l} = x^{2n+1} y x^{-n} y^l x^{-n} y$.

Now we consider surface slope 2.  By the calculations above, $K_{n,0}$ has surface slope 2 with respect to $F$, and $w_{n,0} = x^{2n+1}yx^{-2n}y$.  Applying the automorphism of $\langle x, y \rangle$ which sends $x^{-2n}y$ to $y$ and fixes $x$ sends $w_{n,0}$ to $x^{4n+1} y^2$.  Thus $K_{n,0}$ is $(2,~4n+1)$ Seifert-fibered with respect to $H'$.

Similarly, note that $K_{n,1}$ has surface slope 3, and $w_{n,1} = x^{2n+1} y x^{-n} y x^{-n} y$.  Using the automorphism of $\langle x, y \rangle$ that fixes $x$ and sends $x^{-n} y$ to $y$, we see that $K_{n,1}$ is $(3,~3n+1)$ Seifert-fibered with respect to $H'$.  Thus, we have shown that the slope 2 and slope 3 small Seifert-fibered Dehn surgery on each twist knot can be ``explained" by the primitive/Seifert phenomenon.

\section{Twisted Torus Knots}

\label{s:ttor}
Many examples of primitive/Seifert-fibered knots can be found among the twisted torus knots, which we define in this section.  The simplest twisted torus knots are obtained by adding a full twist to some parallel strands of a torus knot.

Not all twisted torus knots are primitive/Seifert-fibered.  In this section, we develop criteria to determine when a twisted torus knot is primitive or Seifert-fibered with respect to either side of its canonical Heegaard surface.  This involves a detailed analysis of the homotopy class of the twisted torus knot in each handlebody of the Heegaard splitting.

\subsection{Definition and Basic Properties}

\begin{figure}[ht!]
\begin{center}
\includegraphics[width=3in,bb=156 313 439 433]{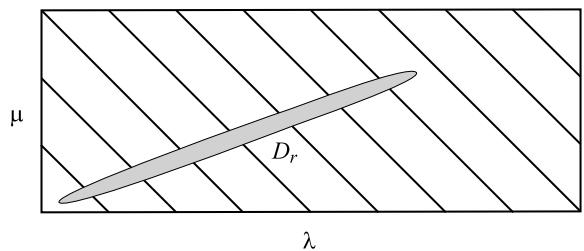}
\caption{The disk $D_r$}
\label{f:baguette}
\end{center}
\end{figure}

Let $\tau$ be an unknotted solid torus contained in $B^3$, the closed $3$-ball, with $\partial \tau$ intersecting $\partial B^3$ in a $2$-disk $D_r$.  Let $\mu$ and $\lambda$ be a meridian-longitude basis for $\partial \tau$.
Let $L$ be a $(p,q)$ torus link contained in $\partial \tau$ that intersects the disk $D_r$ $r$ times as in Figure \ref{f:baguette}, where $0 \leq r \leq p + q$ (in the figure, $p=3$, $q=8$, and $r=7$).  The \emph{torus tangle} $T(p,q)_r$ is formed from $L$ by removing the strands of $L \cap D_r$.  Thus $T(p,q)_r$ is an $r$-string tangle in $B^3$ that is contained in a standard punctured torus (namely $\partial \tau - D_r$) that is properly embedded in $B^3$.

Note that a tangle $T(rp,rq)_r$ with $(p,q) = 1$ is comprised of $r$ parallel strands on the punctured torus $\partial \tau - D_r$, each of which is a $T(p,q)_1$ tangle.  Since $T(rp,rq)_r$ is simply $r$ parallel copies of $T(p,q)_1$, we will denote these torus tangles by $rT(p,q)$.  For convenience, we define $0T(p,q)$ to be the torus tangle with no strings.

A consistent choice of orientations for the components of the original torus link $L$ induces an orientation on the strings of $T(p,q)_r$.  We may form links that are canonically embedded on a genus two Heegaard surface in $S^3$ by gluing together two such torus tangles so that the orientations of the strings match up.  

\begin{dfn}  
\label{d:twtorus}
The \emph{twisted torus knot} $T(p,q) + rT(m,n)$ is obtained by gluing together the tangles $T(p,q)_r$ and $rT(m,n)$ as described above, where $0 \le r \le p+q$, $(p,q)=(m,n)=1$, and $p,q,m \ge 0$.
\end{dfn}

The construction must result in a knot (i.e.\ one component) since the $rT(p,q)$ tangle preserves the gluing pattern of the strings in the original torus knot.

Informally, the twisted torus knots are those knots obtained by splicing together a torus knot and a torus link with $r$ components along $r$ parallel strands on each torus.  We will consider the canonical Heegaard surface to be part of the structure of a twisted torus knot.

\begin{figure}[ht!]
\begin{center}
\includegraphics[bb=202 121 602 467,width=4.5in]{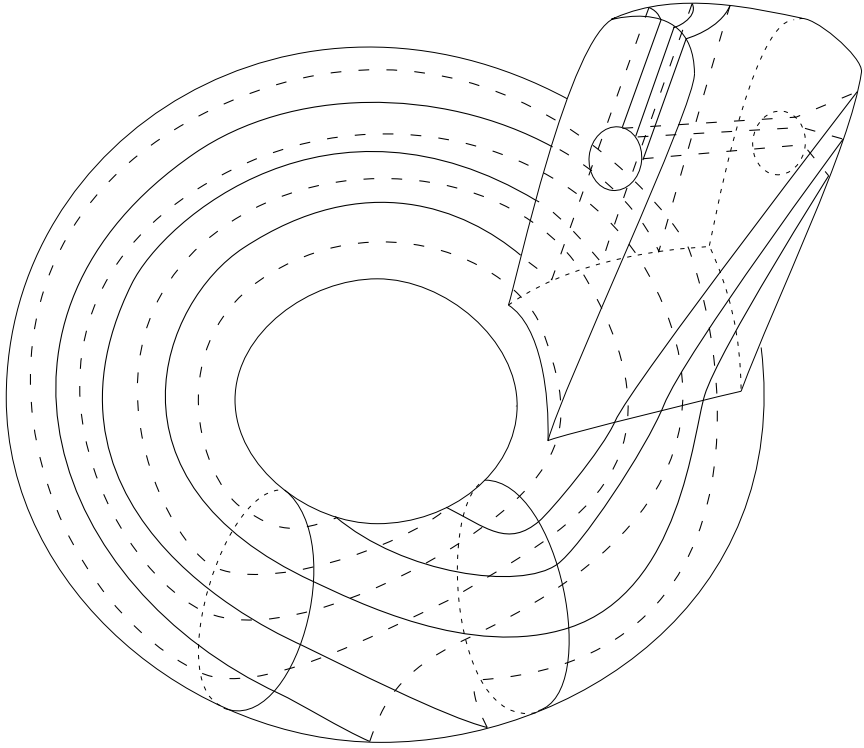}
\caption{The twisted torus knot $T(7,2)+3T(1,1)$}
\label{f:tw_torus}
\end{center}
\end{figure}

Next we calculate the surface slope for a twisted torus knot.  This is the slope at which a small Seifert-fibered surgery may arise.

\begin{pro}
\label{p:tw_torus_slope}
The surface slope of $T(p,q)+rT(m,n)$ with respect to the Heegaard surface described above is $pq +  mnr^2$.
\end{pro}

\proof
Because the ``surface slope'' of an $(a,b)$ torus link is $ab$, each torus tangle $T(a,b)_r$ contributes $ab$ to the surface slope.  The sum of these contributions is $pq + mnr^2$.
\endproof

\subsection{Calculating $w_{p,q,r,m,n}$}

\label{ss:calculatingw}
Now we provide an algorithm to calculate the conjugacy class of a twisted torus knot in the fundamental group of the handlebody $H$.  The algorithm will be used extensively in what follows, and it yields some useful symmetry properties of these conjugacy classes.

Let $w_{p,q,r,m,n}$ and $w'_{p,q,r,m,n}$ be the conjugacy class of the twisted torus knot $K = T(p,q)+rT(m,n)$ in $\pi_1(H)$ and $\pi_1(H')$, respectively.
When any of the values of $p$, $q$, $r$, $m$, or $n$ are apparent or irrelevant, we will omit those subscripts.

\begin{figure}[ht!]
\begin{center}
\includegraphics[bb=202 121 602 467,width=4.5in]{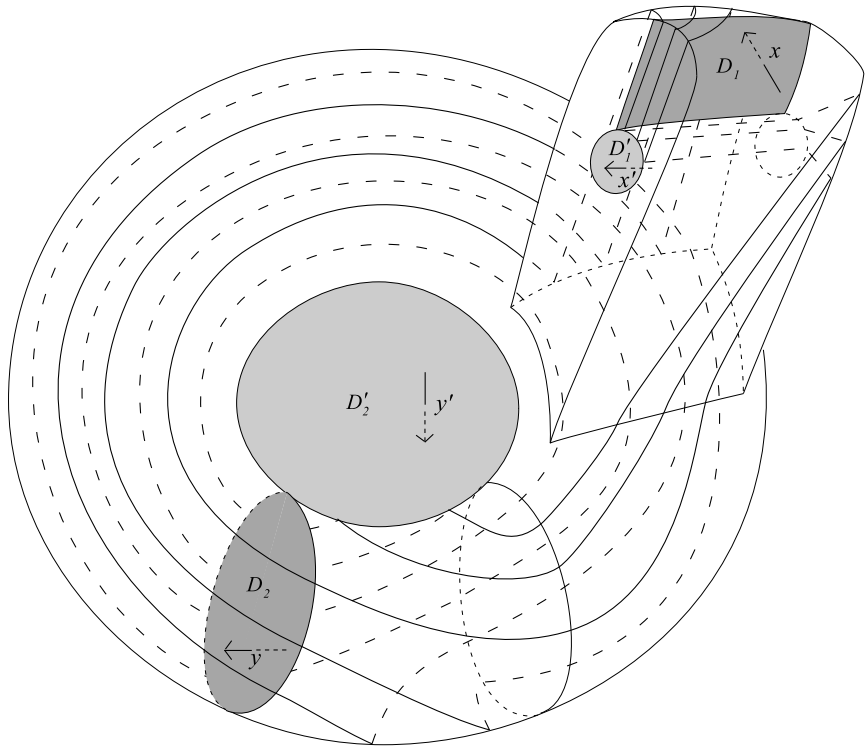}
\caption{Bases for $\pi_1(H)$ and $\pi_1(H')$}
\label{f:bases}
\end{center}
\end{figure}

In the following lemma, we use the disks $D_1$ and $D_2$ in Figure \ref{f:bases} to provide a basis $x$, $y$ for  $\pi_1 (H)$.    Similarly, the disks $D'_1$ and $D'_2$ in Figure \ref{f:bases} describe a basis $x'$, $y'$ for $\pi_1 (H')$.  Note that, with this choice of bases, the exponents of a particular basis element in either $w$ or $w'$ always have the same sign, hence no cancellation can occur.  The disks $D_1$, $D_2$ are distinct from the attaching disk $D_r$ (defined earlier), which in Figure \ref{f:bases} is the rectangular base of the ``mailbox'' that contains $rT(m,n)$.

\begin{lem}
\label{l:wAlgo}
Let $p$, $q$, $r$, $m$, and $n$ be as in definition \ref{d:twtorus}.  Write $r = \bar r + \alpha p$, where $0 \le \bar r < p$ and $\alpha \ge 0$.  The following algorithm calculates $w_{p,q,r,m,n}$ with respect to the basis $x$, $y$ described above.

Mark $p$ points $1,2, \dots p$ consecutively around a circle.  Start at $1$, and jump forward $q$ points.  After each jump, record either $x^{(\alpha + 1)m}y$ or $x^{\alpha m}y$ depending on whether or not the initial point of the jump was between $1$ and $\bar{r}$.  After exactly $p$ jumps, we will have returned to $1$, the starting point. The resulting word in $x$ and $y$ is $w_{p,q,r,m,n}$.
\end{lem}

\begin{figure}[ht!]
\begin{center}
\includegraphics[bb=166 202 451 572,width=3in]{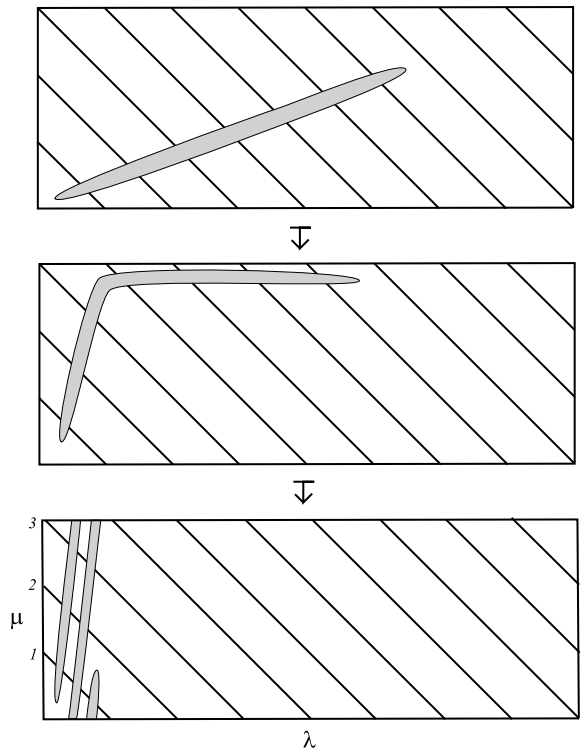}
\caption{Isotoping $D_r$ along the meridinal direction}
\label{f:long_baguette}
\end{center}
\end{figure}

\proof
Let $\tau$ be the solid torus whose boundary contains $T(p,q)_r$ (in Figure \ref{f:bases}, $\tau$ is the large doughnut containing the tangle $T(7,2)_3$).
In order to calculate $w$, it is convenient to isotop $rT(m,n)$ (which is contained in the ``mailbox'' in Figure \ref{f:bases}) so that the attaching disk $D_r$ runs in the meridinal direction of $\partial \tau$.  Figure \ref{f:long_baguette} illustrates how this isotopy can be realized (the figure illustrates the case in which $p=3$, $q=8$, and $r=7$).  Now write $r = \bar{r} + \alpha p$, where $0 \le \bar{r} < p$ and  $\alpha \ge 0$.  After the isotopy, observe that $D_r$ wraps around $\partial \tau$ $\alpha$ complete times in the meridinal direction (intersecting the underlying $(p,q)$ torus knot $\alpha p$ times), then intersecting it $\bar{r}$ additional times.  In Figure \ref{f:long_baguette}, $\alpha = 2$ and $\bar r = 1$.

Now we look at the final frame of Figure $\ref{f:long_baguette}$ and use it to read off $w$.  Start by numbering the intersections of $K$ with $\mu = \partial D_2$ from 1 to $p$ (as in Figure $\ref{f:long_baguette}$).  Then travel along $K$ recording all intersections with $\partial D_1$ as $x$ and $\partial D_2$ as $y$.  Observe that each time $K$ intersects the disk $D_r$ entails one trip along a single strand of the tangle $rT(m,n)$, which results in $m$ consecutive intersections with $D_1$ (hence gives $x^m$).  Note also that each numbered path leaving $\mu$ intersects $D_r$ either $\alpha$ or $\alpha + 1$ times before returning to $\mu$ (contributing either $x^{\alpha m}y$ or $x^{(\alpha + 1)m}y$ to $w$). The numbered paths hit $D_r$ $(\alpha + 1)$ times exactly when the starting point has number less than or equal to $\bar r$. Putting this information together yields the algorithm described in the lemma.
\endproof

\begin{rems}
We have described how to calculate the conjugacy class of the knot with respect to the ``inside'' handlebody $H$.  For the ``outside'' handlebody $H'$, the word $w'_{p,q,r,m,n}$ that $T(p,q)+rT(m,n)$ represents in $H'$ with respect to the basis $x',y'$ pictured in Figure \ref{f:bases} is equal to $w_{q,p,r,n,m}$ with $x$ replaced by $x'$ and $y$ replaced by $y'$. Thus, any result involving $w$ implies an analogous result for $w'$.
\end{rems}

\begin{dfn}  
\label{d:word_equiv}
For $g$ and $h$ in a group $G$, we say $g$ is \emph{equivalent} to $h$ (and write $g \equiv h$) if there is an automorphism of $G$ carrying $g$ to $h$.
\end{dfn}
\begin{lem}
\label{l:symmetries}
The words $w_{p,q,r,m,n}$ enjoy the following properties:
\begin{enumerate}
\item
$w_{p, q, r, m, n}$ does not depend on $n$ (hence we often omit the $n$).
\item
$w_{p,q,r,m} \equiv w_{p,q',r,m}$ if $q \equiv \pm q' \mod p$
\item
$w_{p,q,r,m} \equiv w_{p,q,r',m}$ if $r \equiv \pm r' \mod p$
\end{enumerate}
\end{lem}

\proof
The algorithm described above for calculating  $w_{p,q,r,m,n}$ does not involve $n$ so the first claim is clearly true.  Modifying $q$ by adding a multiple of $p$ obviously has no effect on $w$.  Changing the sign of $q$ results in a conjugate of $w_{p,q,r,m}$, so the second claim is proved.

To prove the third claim in the lemma, recall that the elements in $\pi_1(H)$ that arise in calculating $w$ are $x^{(\alpha + 1)m}y$ and $x^{\alpha m}y$.  The automorphism defined by $x \mapsto x$ and $y \mapsto x^{-\alpha m}y$ carries $x^{(\alpha + 1)m}y$ to $x^my$ and $x^{\alpha m}y$ to $y$.  This shows that $w_{p,q,r,m} \equiv w_{p,q,\bar r,m}$ (as before, $\bar r$ is the remainder of $r$ when divided by $p$).

There is a similar automorphism of $\langle x, y \rangle$ that has the effect of interchanging $x^{\alpha m}y$ and $x^{(\alpha + 1)m}y$,
namely $x \mapsto x^{-1}$ and $y \mapsto x^{(2\alpha + 1)m}y$
This shows that $w_{p,q,r,m} \equiv w_{p,q,p-\bar r,m}$, so the lemma is proved.
\endproof

\subsection{Primitivity of Twisted Torus Knots}

We develop a simple criterion for determining whether or not a twisted torus knot is primitive with respect to either handlebody of its canonical Heegaard splitting.  Since this only depends on $w_{p,q,r,m,n}$, we focus on these elements of the free group on two generators.

We will use the following necessary condition for primitivity in the proof of the next theorem. Up to cyclic reordering and the automorphism that interchanges $x$ and $y$, a primitive word in the free group on $x$ and $y$ with positive exponents has one of the following regular forms \cite{cmz:1981}:

\begin{itemize}
\item
$x^ly$,
\item
$x^{l_1}yx^{l_2}y \dots x^{l_k}y$ where $\{l_i\} = \{e, e+1\}$ for some positive integer $e$.
\end{itemize}

\begin{rems}
The converse is not true.  For example, $xyxyx^2yx^2y$ is in regular form, but the fact that it is equivalent to $x^2y^2$ under an automorphism of $\langle x, y \rangle$ shows that it is not a primitive element.
\end{rems}

\begin{thm}
\label{t:prim}
$w_{p,q,r,m}$ is primitive if and only if 
\begin{enumerate}
\item
$p=1$; or
\item
$m=1$ and $r \equiv \pm 1$ or $\pm q \mod p$.
\end{enumerate}
\end{thm}

\proof
If $p=1$ then $w$ is obviously primitive since it has the form $x^ly$.  From now on we assume $p>1$.

If $r=0$ then $w$ is a primitive element raised to the $p$th power.
Since $p\ge 2$, clearly $w$ is not primitive.  Henceforth we assume that $r>0$.

Now we show that if $p>1$, $r>0$, and $m>1$ then $w$ is not primitive.  If $\alpha = 0$, then the fact that $r \ne 0$ and $m>1$ implies that $w$ contains powers of both $x$ and $y$ that are greater than $1$, so $w$ has no regular form.  If $\alpha > 0$ then the exponents of $x$ differ by $m$, which is greater than $1$, so $w$ has no regular form.

What remains is the main case---when $p>1$, $r>0$, and $m=1$.
By lemma \ref{l:symmetries} we may also assume that $r < p$ and $1 < q < p/2$.  Thus $\alpha = 0$, so our ``building blocks'' $x^{(\alpha +1)m}y$ and $x^{\alpha m}y$ simplify to $xy$ and $y$.  For the remainder of this proof we apply to $w$ the automorphism $xy \mapsto x$ and $y \mapsto y$ of $\pi_1(H) = \langle x, y \rangle$.  We will divide the rest of the proof into the following four cases.

\rk{\bf Case 1} $1 < r < q$

\begin{clm}
For $1 < r < q$, the exponent of each power of $x$ in $w$ is one.  However, the set of exponents of $y$ that appear in $w$ contains integers that differ by more than one, so $w$ has no regular form, hence cannot be primitive.
\end{clm}

\noindent
Recall our algorithm for calculating $w_{p,q,r,m}$.  Since $1 < r < q < p/2$, each trip around the circle encounters at most one $x$ and at least one $y$. This implies that $w$ has the form $xy^{a_1}  \dots  xy^{a_r}$ where each $a_i$ is a positive integer.

We define $\sigma$ to be the smallest number of jumps required for any element of $\{1, \dots , r\}$ to return to $\{1, \dots , r\} $.  That is
$$\sigma = \min_{a \in \{1, \dots , r\} } \min_{j \in \mathbb{N}} \{ j : (a + jq) \bmod p \in \{1, \dots , r\} \} $$

By the comments above, $\sigma \ge 2$.  In addition, $\sigma-1$ is the smallest exponent of $y$ that appears in $w$ as a cyclic word.  In particular, up to cyclic reordering, $w$ has a subsequence of the form $xy^{\sigma-1}x$.

Choose $\Sigma \in \{1,  \dots , r\}$ that realizes the initial element for $\sigma$, and let $\Theta$ be the element in $\{1,  \dots , r\}$ to which it returns first (thus $\Theta \equiv \Sigma + \sigma q \mod p$). $\Sigma \ne \Theta$ since $r>1$.

Suppose $\Theta < \Sigma$.  Consider the following sequence $\bmod p$
$$\Sigma-\Theta \mapsto \Sigma-\Theta +q \mapsto \dots \mapsto \Sigma-\Theta +\sigma q \mapsto \Sigma-\Theta +(\sigma+1)q.$$
By construction, the first term is in $\{1,  \dots , r\}$.  By the minimality of $\sigma$, the next $\sigma-1$ terms cannot be in $\{1,  \dots , r\}$.  The next term,  $\Sigma-\Theta +\sigma q$, equals $0$ modulo $p$ by the definition of $\Theta$, and the last term is not in $\{1,  \dots , r\}$ since $q>r$ and $q<p/2$. Thus we have a subsequence in $w$ of the form $xy^{\sigma+1}$, one of the form $xy^{\sigma-1}x$, and all exponents in $w$ are positive.  Thus, when $\Theta < \Sigma$, we have shown that $w$ satisfies the claim, and cannot be primitive.

When $\Theta > \Sigma$, we consider a sequence similar to the one above that begins with the term $(r + 1) - (\Theta - \Sigma)$.  The same type of analysis again shows that there is a subsequence of the form $xy^{\sigma+1}$, so we conclude that $w$ is not primitive.

\rk{\bf Case 2} $q < r < p - q$

\begin{clm} 
When $q < r < p - q$, $w$ is a word in $\langle x, y \rangle$ with positive exponents such that there is an exponent of both $x$ and $y$ that is greater than one.  Thus no combination of conjugation and interchanging $x$ and $y$ will put $w$ in regular form, so $w$ can not be primitive.
\end{clm}

\noindent
Since $q$ is less than $r$, at least two $x$ terms in a row are obtained on the first trip around the circle.  Similarly, since $q < p - r$, when the point lands on $r + 1$, the next point hit is less than or equal to $p$, so at least two $y$ terms are obtained.

\rk{\bf Case 3} $p-q < r < p-1$

\noindent
$w$ is not primitive by case 1 and the symmetry described in Lemma \ref{l:symmetries}.

\rk{\bf Case 4} $r = 1$, $q$, $p-q$, or $p-1$

\noindent
When $r = 1$, $w=xy^k$ for some integer $k$, so $w$ is primitive.  When $r=q$, $w_{p,q,q,1}$ is exactly the word in $\langle x, y \rangle$ described in \cite{oz:1981}, where it is proved to be the unique primitive element (up to conjugacy) in $\langle x, y \rangle$ with abelianization $(r,p-r)$.  By Lemma \ref{l:symmetries}, $w$ is also primitive when $r=p-1$ or $r=p-q$.
\endproof

\subsection{Which $w_{p,q,r,m}$ are Seifert-fibered?}

\label{ss:sfwords}
Here we give criteria which allow one to determine to a great extent which $w_{p,q,r,m}$ are Seifert-fibered.  In particular, Propositions \ref{p:sf_hyper}, \ref{p:middle_sf}, and \ref{p:sf_end} describe three circumstances in which $w$ is Seifert-fibered.  These three types of Seifert-fibered elements are summarized in Table\ref{tb:sfwords} at the end of this section.  Moreover, we conjecture that the three types describe all $w_{p,q,r,m}$ that are Seifert-fibered.

First we state some results concerning presentations of the groups 
$$G_{a,b} = \langle x, y \ | \ x^ay^b \rangle.$$
Consider two sets of  conjugacy classes of elements $w_1, \dots w_k$ and $w_1', \dots, w_k'$ in a free group $F$.  The two sets of conjugacy classes are Nielsen equivalent if $F$ has an automorphism mapping each $w_i$ to a conjugate of $w_i'$ (where the conjugating factor depends on $i$).  Two presentations $\langle S_1, \dots, S_n \ | \ R_1(S), \dots R_l(S) \rangle$ and $\langle S_1, \dots, S_n \ | \ R_1'(S), \dots R_l'(S) \rangle$ are Nielsen equivalent if the set of conjugacy classes $R_1(S), \dots R_l(S)$ is Nielsen equivalent to $R_1'(S)^{\pm 1}, \dots R_l'(S)^{\pm 1}$ in the free group $\langle S_1, \dots, S_n \rangle$.

Work of Zieschang and Collins completely determined the Nielsen equivalence classes of 1-relator presentations of the groups $G_{a,b}$ \cite{zieschang:1977, collins:1978}.  Let $v_{s,t}(x,y)$ be the unique primitive element up to conjugacy in $\langle x, y  \rangle$ with $(s,t)$ as its abelianization (see above and \cite{oz:1981}).  If $u$ and $u'$ are elements of $\langle x, y \rangle$, then let $v_{s,t}(u, u')$ denote the conjugacy class in $\langle x,y \rangle$ obtained by substituting $u$ for $x$ and $u'$ for $y$ in $v_{s,t}(x,y)$.  Then the following are the Nielsen equivalence classes of one-relator presentations of the groups $G_{a,b}$:
$$ \langle x, y | v_{a,k}(x, y^b) \rangle $$
$$ \langle x, y | v_{l,b}(x^a, y) \rangle $$
where $(k,a)=(l,b)=1$, $0 < 2k < a$, and $0 < 2l < b$.
Note that larger values of $k$ or $l$ still yield a valid (but redundant) presentation of $G_{a,b}$.  This provides a characterization of the Seifert-fibered words in $\langle x, y \rangle$.

First we need a lemma which shows that the elements $w_{p,q,r,1}$ are the prototypes for all the $w_{p,q,r,m}$.  The lemma follows immediately from the algorithm already given for calculating $w_{p,q,r,m,n}$.  For the rest of this section we may assume, as usual, that $1 \le q < p/2$ and $1 \le r \le p$.  

\begin{lem}
\label{l:endo}
Let $\mathbf{a_m}$ be the endomorphism of $\langle x,y \rangle$ given by $x \mapsto x^m$ and $y \mapsto y$.  Then $w_{p,q,r,m} =  \mathbf{a_m}(w_{p,q,r,1})$.
\end{lem}

\begin{pro}
\label{p:sf_hyper}
If $w_{p,q,r,1}$ is primitive then $w_{p,q,r,m}$ is $(m, p)$ Seifert-fibered.
\end{pro}

\proof
In our original choice of basis for $\pi_1(H)$, the abelianization of $w_{p,q,r,1}$ is $(r,p)$.  By lemma \ref{l:endo}, and the characterization Seifert-fibered words above it is immediate that $w_{p,q,r,m} = \mathbf{a_m} (w_{p,q,r,1})$ is $(m,p)$ Seifert fibered.
\endproof

\begin{cor}
All $1$-relator presentations of the groups $G_{a,b}$ can be realized geometrically, \emph{i.e.} as the obvious $1$-relator presentation of $\pi_1$ for a $3$-manifold obtained by adding a $2$-handle to a handlebody of genus $2$.  
\end{cor}

This Corollary appears (implicitly) in \cite{brz:1988}.

\begin{pro}
\label{p:middle_sf}
For any integer $\beta$ such that $1 \le \beta < p/q$, $w_{p,q,\beta q,1}$ and $w_{p,q,p-\beta q,1}$ are $(\beta, p - \beta q)$ Seifert-fibered.
\end{pro}

For the rest of this section, in describing $w_{p,q,r,1}$ we will use the basis for $\pi_1(H)$ described before Case 1 in the proof of Theorem \ref{t:prim}. 

\begin{figure}[ht!]
\begin{center}
\includegraphics[bb=111 414 549 451,width=5in]{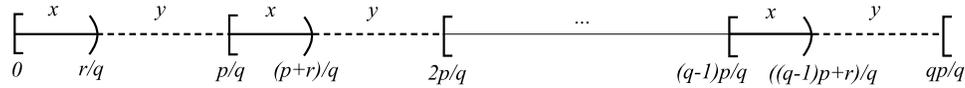}
\caption{Another description of $w_{p,q,r,1}$}
\label{f:numline}
\end{center}
\end{figure}

Before we begin the proof of Proposition \ref{p:middle_sf}, we give an alternate description of the words $w_{p,q,r,1}$.   We may divide up the interval $[0,p)$ in $\mathbb{R}$ into $q$ intervals $[ip/q, (ip+r)/q)$, and $q$ intervals $[(ip+r)/q, (i+1)p/q)$, where $0 \le i < q-1$ in each case (as in Figure \ref{f:numline}).  Note that every integer from $0$ to $p-1$ falls into a subinterval of the first type or of the second type.  We label each integer that falls into one of the $q$ subintervals in the first list by $x$ and those that fall into one of the $q$ subintervals of the second type by $y$ (see Figure \ref{f:numline}).  Note that, depending on $r$, $p$, and $q$, a given subinterval might contain several integers (or none).  We claim that reading off the sequence of $x$'s and $y$'s associated to the integers $0, 1, \dots, p-1$ gives the word $w_{p,q,r,1}$.  We will use this description in the proof of the next proposition.

To see why this agrees with the description we gave for $w_{p,q,r,1}$ in Lemma \ref{l:wAlgo}, consider the covering map $\rho_1: \mathbb{R} \rightarrow S^1$, given by identifying each real number $u$ with $u+p$, and the $q$-fold covering map $\rho_2: S^1 \rightarrow S^1$ given by $u \mapsto qu$.  Now consider the image of the points $0, 1, \dots, p-1$ under the composite map $\rho_2 \circ \rho_1:\mathbb{R} \rightarrow S^1$.  Replace each integer $j$ in the sequence $0, 1, \dots, p-1$ with  $x$ or $y$ depending on whether or not $0 \le (\rho_2 \circ \rho_1) (j) < r$ (here we use the coordinates $[0,p)$ for $S_1$, obtained from $\rho_1$).  Note that this condition for replacing each integer in $0, \dots, p-1$ with $x$ or $y$ agrees exactly with that described in Lemma \ref{l:wAlgo}, because we may think of $j$ as representing the number of jumps taken (and the mapping $\rho_2$ ``executes'' the jumps).  Looking at $(\rho_2 \circ \rho_1)^{-1}([0,r))$ gives the subintervals of the of the  first type described above, while $(\rho_2 \circ \rho_1)^{-1}([r,p))$ gives the subintervals of the second type.  Translating the condition $0 \le (\rho_2 \circ \rho_1) (j) < r$ to one concerning the (lifted) subintervals gives the description of $w_{p,q,r,1}$ from the previous paragraph.

We will also need the following easy lemma:

\begin{lem}
\label{l:aut}
Suppose $u=xy^{k_1} \dots xy^{k_j}$ is an element in $\langle x, y \rangle$. Let $\tilde u = xy^{k_1-l} \dots xy^{k_j-l}$.  Then $u \equiv \tilde u$; in particular, $u$ is primitive if and only if $\tilde u$ is.
\end{lem}

\proof[Proof of Proposition \ref{p:middle_sf}]
In the case $r=q$ we know from the last section that $w_{p,q,q,1} = xy^{k_1} \dots xy^{k_q}$ is the unique primitive word (up to conjugacy) in $\langle x, y \rangle$ with abelianization $(q, p-q)$.  Each $k_i$ is greater than zero. Let $\beta$ be an integer between $1$ and $p/q$. Using the new description of the words $w_{p,q,r,1}$ given above, it is easy to see that changing $r$ from $q$ to $\beta q$ changes each exponent of $x$ from 1 to $\beta$, and subtracts $\beta - 1$ from each exponent of $y$.  So if $w_{p,q,q,1}$ is as above, we have:
$$w_{p,q,\beta q,1}  = x^{\beta}y^{k_1+1-\beta} \dots x^{\beta}y^{k_q+1-\beta}.$$
Note that the abelianization of $w_{p,q,\beta q,1}$ is $(\beta q, p-\beta q)$.

By the previous lemma, and the observations above, $w_{p,q,\beta q,1}$ is equivalent to the $(\beta, p-\beta q)$ Seifert-fibered word $v_{q,p-\beta q}(x^\beta, y)$.

By symmetry, $w_{p,q,p-\beta q,1}$ is also $(\beta, p-\beta q)$ Seifert-fibered.
\endproof

Next we describe a third type of $w_{p,q,r,m}$ which is Seifert-fibered. Let $\hat q^{-1}$ be the smallest positive integer congruent to $\pm q^{-1}$ modulo $p$.
Then we have the following result, where $\{\ \}$ denotes the least integer function.

\begin{pro}
\label{p:sf_end}
For $1 \le r \le \{p/ \hat q^{-1}\}$, $w_{p,q,r,1}$ is $(r, p-r\hat q^{-1})$ Seifert-fibered.  By symmetry, so is $w_{p,q,p-r,1}$. 
\end{pro}

\proof
Recall the first method we described for computing $w_{p,q,r,1}$ involving the integers from 0 to $p-1$ arranged in a circle. We know that we start at 0, and return to 0 only after $p$ jumps of length $q$.  Each time we land on a point between 0 and $r-1$, we record an $x$, while the other points result in a $y$.  We may write $w_{p,q,r,1} = xy^{k_0} \dots xy^{k_r-1}$, where the $k_i$ are non-negative integers.  The exponents of $y$ are the number of jumps required to return to the interval $[0, r-1]$ during the orbit.

Consider the following congruences:
\begin{align*} 
     j_0 q & \equiv 0 \mod p\\
     j_1 q & \equiv 1 \mod p\\
             &\vdots\\
     j_{r-1} q &\equiv r-1 \mod p
\end{align*}
So $j_i \equiv i q^{-1} \bmod p$.  If we take $j_i$ to be the smallest positive solution to the congruence, then $j_i$ is the minimum number of jumps of length $q$ required to land on the point $i$, having started at 0.

Sort the $j_i$ in monotonic order and reindex (if necessary) to reflect the new ordering.  This sorted list describes the total number of jumps that have occurred each time that the point lands in the interval $[0, r-1]$ during the orbit, having started at zero.  The differences between consecutive $j_i$ give the exponents of $y$ in $w_{p,q,r,1}$, in particular:
$$k_i = j_{i+1}-j_i-1$$
where we consider the subscripts to be in $\mathbb{Z}_r$.

The $j_i$ are \emph{already} in monotonic order after being reduced modulo $p$ exactly when
$1 \le r \le \{p/ \hat q^{-1}\}$.  Thus  $k_i = \hat q^{-1}-1$ for $i = 0, \dots, r-2$ and 
$$k_{r-1} = p - (r-1) \hat q^{-1}-1.$$
Hence 
$$ w_{p,q,r,1} = (xy^{\hat q^{-1}-1})^{r-1}xy^{p-(r-1) \hat q^{-1}-1}.$$
If we apply the automorphism defined by $xy^{\hat q^{-1}-1} \mapsto x$ and $y \mapsto y$, we find that
$$ w_{p,q,r,1} \equiv x^ry^{p-r \hat q^{-1}}. $$
We conclude that, for $1 \le r \le \{p/ \hat q^{-1}\}$, $w_{p,q,r,1}$ is $(r, p-r \hat q^{-1})$ Seifert-fibered.
\endproof

Table \ref{tb:sfwords} summarizes the three types of Seifert-fibered $w_{p,q,r,m}$ described in this section.  We call the first type of Seifert-fibered $w_{p,q,r,m}$ \emph{hyper Seifert-fibered} since they arise from a primitive twisted torus knot with one full twist (i.e.\ $m=1$) by simply increasing the number of full twists.  The remaining types have only one full twist.  Since $m=1$ for these types, we fix $p$ and $q$, and consider the integers $\bmod p$ lined up from $1$ to $p$.  We must choose $r \bmod p$ from these to make $w_{p,q,r,m,}$ Seifert-fibered.  The values of $r$ are evenly spaced throughout this range for the second type, so they are called \emph{middle Seifert-fibered}.  The values of $r$ are clustered at both ends for the third type so they are called \emph{end Seifert-fibered}.

\begin{table}
\begin{center}
\begin{tabular} {|l||l|l|l|} \hline
\strut type   & critical fibers      &  $w_{p,q,r,m}$ satisfying\\ \hline \hline
\strut hyper  & $(p,m)$              &  $|m| > 1$, $r \equiv \pm 1$ or $\pm q \bmod p$               \\ \hline
\strut middle & $(\beta, p-\beta q)$ & $m=1$, $r\equiv \pm \beta q \bmod p$, where $1 \le \beta < p/q$ \\ \hline
\strut end    & $(r,p-r\hat q^{-1})$ & $m=1$, $r\equiv \pm \bar r \bmod p$, where $1 \le \bar r \le \{p/\hat q^{-1}\}$  \\ \hline
\end{tabular}
\caption{The three types of Seifert-fibered elements}
\label{tb:sfwords}
\end{center}
\end{table}

\section{Primitive/middle-SF Knots and Surgeries}

\label{s:class&mult}
The results of the last section provide much information about which twisted torus knots are primitive or Seifert-fibered with respect to either handlebody of their associated Heegaard splittings.  The first part of this section is devoted to identifying a large class of twisted torus knots that are \emph{simultaneously} primitive on one side and Seifert-fibered on the other.

In the previous section,  the surgery slope and the multiplicities of two of the three critical fibers in the surgered manifold were calculated.  In the last part of this section we determine the multiplicity of the third critical fiber.  To do this, we will need to find a curve in $\partial H - N(K)$ which becomes an ordinary fiber in $H \cup_K \text{2-handle}$ (the latter is a Seifert-fibered space over $D^2$ with two critical fibers).

Many twisted torus knots are doubly Seifert-fibered.  As mentioned in section \ref{ss:psf}, the latter must have Dehn surgeries that are graph manifolds or Seifert-fibered spaces with base $S^2$ and four critical fibers.  An analysis like that in this chapter would determine which graph manifolds (and possibly Seifert-fibered spaces) arise from Dehn surgery on doubly Seifert-fibered twisted torus knots.

\subsection{Some primitive/Seifert-fibered twisted torus knots}

We use the criteria for recognizing primitive or Seifert-fibered twisted torus knots developed in the last section to find many twisted torus knots which are primitive/Seifert-fibered.  

Our goal is to show that primitive/Seifert-fibered knots are abundant in some sense.  To do this, it is enough to focus on twisted torus knots formed using a single full twist (i.e.\ $K = T(p,q) + rT(1,\pm 1)$), and we assume that $r < \max \{p, q\}$.  Within this class of twisted torus knots, we determine exactly which ones are primitive/\emph{middle}-Seifert-fibered.  A similar analysis could be carried out for primitive/end-Seifert-fibered and primitive/hyper-Seifert-fibered twisted torus knots.

For convenience, we will write $K(p,q,r,m,n)$ in place of $T(p,q) + rT(m,n)$.  In particular, we will be focusing on the knots $K(p,q,r,1,\epsilon)$, where $\epsilon = \pm 1$.

Without loss of generality, we may assume that $K$ is Seifert-fibered with respect to $H$ and primitive with respect to $H'$.  We will not include the doubly primitive case, or any twisted torus knot that is obviously a torus knot (i.e.\ $r=1$, $p$, or $q$; or $q=1$).

\begin{thm}
\label{t:class}
The twisted torus knots $K(p,q,r,1,\epsilon)$ with $r < \max \{p, q\}$ and $\epsilon=\pm 1$ that are middle Seifert-fibered with respect to $H$ and primitive with respect to $H'$ are given by the following values of $(p,q,r)$:
\begin{center}
\begin{tabular} {|l||l|l|} \hline
\strut   & $(p,q,r)$     & {\rm satisfying}\\ \hline \hline
\strut 1 &  $(p,q,2q-p)$ & $\frac{p+1}{2} < q < p$\\ \hline
\strut 2 &  $(p,q,p-kq)$ & $1 < q < \frac{p}{2}$, $2 \le k \le \frac{p-2}{q}$\\ \hline
\strut 3 &  $(ls+l+\delta, ls+\delta, ls)$ & $s \ge 2$, $l \ge 2-\delta$, $\delta=\pm 1$\\ \hline
\strut 4 &  $(p,tp-l,tp-l-1)$, where $p=ls+l+1$ & $s \ge 2$, $l \ge 1$, $t \ge 2$\\ \hline
\strut 5 &  $(p, s+tp, s-1+tp)$, where $p=ls-1$ & $s \ge 3$, $l \ge 3$, $t \ge 1$\\ \hline
\end{tabular}
\end{center}
\end{thm}

\begin{rems}
Using the results of the previous section, it is easy to verify that each of the knots in the theorem is indeed primitive/middle-Seifert-fibered.
\end{rems}

\proof
As usual, let $\bar{q}$  ($\bar r$ respectively) be the smallest positive integer congruent to $q$ modulo $p$ ($r$ modulo $p$, respectively).  Let $\hat q$ be the smallest positive integer congruent to $\pm q$ modulo $p$.

\rk{\bf Case 1} $r \equiv \pm p \mod q$

\noindent
We show that the knots in 1 and 2 are the only possibilities.

If $q > p$, then $r < q$. Thus $r \equiv \pm p \mod q$ implies that $r=p$ or $r=q-p$.  The first we ignore since it describes a torus knot, and the second because it describes a doubly primitive knot with respect to the given Heegaard surface.

We consider the case $q < p$.  Since $K$ is middle Seifert-fibered, $r \equiv \pm kq \mod p$ for some integer $k$ satisfying $1 < k < p/q$.  It is not hard to show that the integers that satisfy these conditions and $r \equiv \pm p \mod q$ have the form $\pm kq \pm p$.  Since we require $r < \max \{p,q\}$ the only possibilities are $\pm (kq -p)$.

Now if $q > p/2$, then $r= \pm (kq - p)$ and $p > r > 0$ imply that $r=p-q$ or $r=kq-p$ where $k = 2$ or 3.  The former implies that $K$ is doubly primitive, so we discard that solution.  Since $p> q > p/2$ we have that $\hat q = p-q$.  Hence $k=2$ gives a legitimate solution because 
$$r = 2q - p = 2(p- \hat q) - p = p - 2 \hat q.$$
In this case $K$ is $(2, 2q -p)$ Seifert-fibered (with respect to $H$).  Since we discard any torus knot solutions, we require that $q > (p+1)/2$.  These are the first class of knots described in the theorem.

If $k=3$ we have
$$ r= 3q-p = 3(p- \hat q) - p = 2p - 3\hat q.$$
For $K$ to be middle Seifert-fibered, we must have that either $r = a\hat q$ or $r=p-a\hat q$ for some integer $2 \le a \le p/\hat q$.  This implies that either $2p = (a+3)\hat q$ or $p=(3-a)\hat q$.  Since $(p, \hat q) =1$ the only possibilities are $\hat q =1$ or 2 in the first case and $\hat q =1$  in the second case.  Since $\hat q = 1$ or 2, there are only a few values of $p$ for which $r = 2p -3\hat q$ takes on values between 1 and $p-1$.  Only one actual solution arises and it is doubly primitive.

Now we consider $q < p/2$.  As before, we have the possibilities $r = \pm (kq-p)$. If $r = p - kq$ and  $2 \le k \le (p-2)/q$ then we get the second set of solutions described in the theorem.  We show that no other values of $r$ give solutions.  The remaining case is when $r=kq-p$ and $r=j\hat q$ for some integer $j$.  These equations imply that $\hat q | p$. Since $(\hat q,p)=1$ this means that $\hat q$ must equal 1.  Because $q < p/2$, $\hat q = q$ so we conclude that $q = 1$, hence the knot is a torus knot.

\rk{\bf Case 2} $r \equiv \pm 1 \mod q$

\noindent
We have described all primitive/middle-Seifert-fibered knots where the primitive side satisfies $r \equiv \pm p \mod q$.  Now we consider those whose primitive side satisfies $r \equiv \pm 1 \mod q$.  When we find solutions for which $p \equiv \pm 1 \mod q$, they will necessarily have been described already since this overlaps with the case $r \equiv \pm p \mod q$ already studied.

The following lemma explains the structure of those solutions in case 2 with $q>p$.

\begin{lem}
\label{l:fund_soln}
$K(p,q,r,1,\epsilon )$ is primitive/middle-Seifert-fibered with $q>p$ and $r \equiv \pm 1 \mod q$ if and only if $r=1$ or $q-1$ and $K(p,\bar q, \bar r, 1, \epsilon )$ is primitive/mid\-d\-le-Seifert-fibered.
\end{lem}

\proof[Proof of Lemma \ref{l:fund_soln}]
Suppose that $K(p,q,r,1,\epsilon )$ is primitive/middle-Seif\-ert-fibered.  Since $q>p$, we know that $1 \le r < q$.  Thus $r \equiv \pm 1 \mod q$ implies that $r=1$ or $q-1$.  

By the symmetries in Lemma \ref{l:symmetries}, $w_{p, \bar q, \bar r,1} \equiv w_{p,q,r,1}$, so $K(p,\bar q, \bar r,1,\epsilon )$ is middle Seifert-fibered with respect to $H$.  By Theorem \ref{t:prim}, $w_{\bar q, p, \bar r, \epsilon }$ is primitive since $\bar r = 1$, or $\bar q - 1$.

The reverse implication is similar. 
\endproof

Since $r=1$ corresponds to a torus knot, we will only consider solutions arising from the lemma with $r=q-1$.

This lemma simplifies the remainder of our task.  First we find all solutions with $q<p$, keeping in mind those which have the property that $r=q-1$.  Then, to generate all solutions with $q>p$, we need only add a fixed multiple of $p$ to both $q$ and $r$ for each such solution found with $q<p$.

Now we assume that $q<p$.  Since $1 \le r < p$ and $r \equiv \pm 1 \mod q$, the candidates for $r$ are $jq \pm 1$ where $j \ge 0$ and $jq \pm 1 < p$.

\rk{\bf Subcase i}  $q < p/2$.

\noindent
Since $\hat q = q$, and $r=jq \pm 1$ as above, we have that either $jq \pm 1 = sq$ or $jq \pm 1 = p-sq$ for some integer $1< s < p/q$.

The first equation implies that $q=1$ in which case the knot is a torus knot.  Note that this solution does not generate any solutions with $q>p$ as in the previous lemma.

The second equation implies that $p \equiv \pm 1 \mod q$. But then 
$$r \equiv \pm p \equiv \pm 1 \mod q.$$
Thus these solutions fall in the overlap between Case 1 and Case 2.  In fact, they are included in the second type of solution described in the theorem.  Thus they contribute no new solutions when $q<p$, but we will need to consider the solutions that they generate via the previous lemma. 

 We will return to this later.

\rk{\bf Subcase ii} $q>p/2$.

\noindent
The only possible values of $r=jq \pm 1$ that lie between 1 and $p-1$ are $r=1$, $q-1$ or $q+1$.  The first possibility describes a torus knot, and the solutions with $q>p$ obtained from this using the previous lemma are all doubly primitive, so we discard them. The other possibilities are
\begin{align*}
q \pm 1 & = p - s\hat q\\
\text{and } q\pm 1  &= s\hat q,
\end{align*}
where $1 \le s < p/\hat q$.  Substituting $\hat q = p - q$ yields the following:
\begin{align*}
\pm 1 & = (p-q)(1-s)\\
\text{and } \pm 1  &= sp-(1+s)q.
\end{align*}

The only solution to the first equation that satisfies $q<p$ and $s>0$ is
$$K(p, p-1, p-2, 1,\epsilon ).$$  The solutions that they generate via the lemma are subsumed by the fourth class of knots in the theorem, which we will discuss shortly.

It is easy to show that the following are all integer solutions to the second equation:
$$ (p,q) = \mp (1+l(s+1), 1+ls),$$
where $l \in \mathbb{Z}$.  Taking into account that $p$, $q$, and $s$ are required to be positive (and replacing $l$ by $-l$ in the ``$-$" case), these solutions describe are the knots
$$K(ls+l+\delta, ls+\delta, ls,1,\epsilon).$$
where $l, s > 0$, and $\delta = \pm 1$.  These are exactly the third category of knots described in the theorem.  The additional restrictions on $l$ and $s$ in the theorem ensure that $K$ is not doubly primitive.  Note that when $l=\delta=1$, we get exactly the solutions described in the previous paragraph.

The solutions in the third category induce solutions with $q>p$ exactly in the case when $\delta = 1$.  These are the solutions in the fourth category in the theorem.

Thus far we have described every category of knots in the theorem except for the fifth, which are obtained by applying the lemma to all those solutions in category 2 which satisfy $p \equiv \pm 1 \mod q$ and $r=q-1$.
\endproof

\subsection{Locating the Ordinary Fiber}

We already know the homology and the multiplicities of two critical fibers for the small Seifert-fibered surgeries arising from the knots in Theorem \ref{t:class}.  A first step towards determining the multiplicity of the third critical fiber is to find a curve in the boundary of the relevant handlebody that is disjoint from $K$ and becomes an ordinary fiber once a $2$-handle is attached along $K$.

In fact, the only information that we will need to extract from this curve is its homology class in the complementary handlebody.

We assume in this subsection that $K=K(p,q,r,1,\epsilon)$ is middle-Seifert-fibered with respect to $H$.  We will describe explicitly a curve $f$ in the boundary of $H$ that is disjoint from $K$ and is isotopic to an ordinary fiber in
$$W = H \cup_K \text{2-handle}.$$

\noindent
\rk{Notation}  Write $q=\tilde q + \gamma p$ where $-p/2 < \tilde q < p/2$ and $\gamma \ge 0$.  Let $\hat q = |\tilde q|$.  We may write $r = \bar r + \beta p$ where $0< \bar r < p$ and $\beta \ge 0$.

The curve $f$ that we describe will be somewhat different depending on whether or not $\tilde q < 0$ and whether $r=\alpha \hat q$ or $r= p - \alpha \hat q$ where $1 < \alpha < p/q$.  This gives four natural cases.  First we calculate the homotopy class in $\pi_1(W)$ of an ordinary fiber.  Then we describe the curve $f$ explicitly for one example in each case.  Since the dynamics of the curve on the genus two surface are specified by which of the four cases we are in, one example essentially describes the general case.  It is then a simple matter to calculate the homology class that $f$ represents in $H'$ with respect to our usual basis.

When $\bar r = p -\alpha \hat q$, then $w =v_{\hat q, p-\alpha \hat q}((x^{\beta}y)^{\alpha}, x^{\beta+1}y)$ and the ordinary fiber is $(x^{\beta}y)^{\alpha}$ in $\pi_1(W)$, where the generators are induced from those for $H$.  When $\bar r = \alpha \hat q$ then we need only interchange the roles of $\beta$ and $\beta + 1$ in both expressions.

\begin{figure}[ht!]
\begin{center}
\includegraphics[bb=98 194 501 633,width=4in]{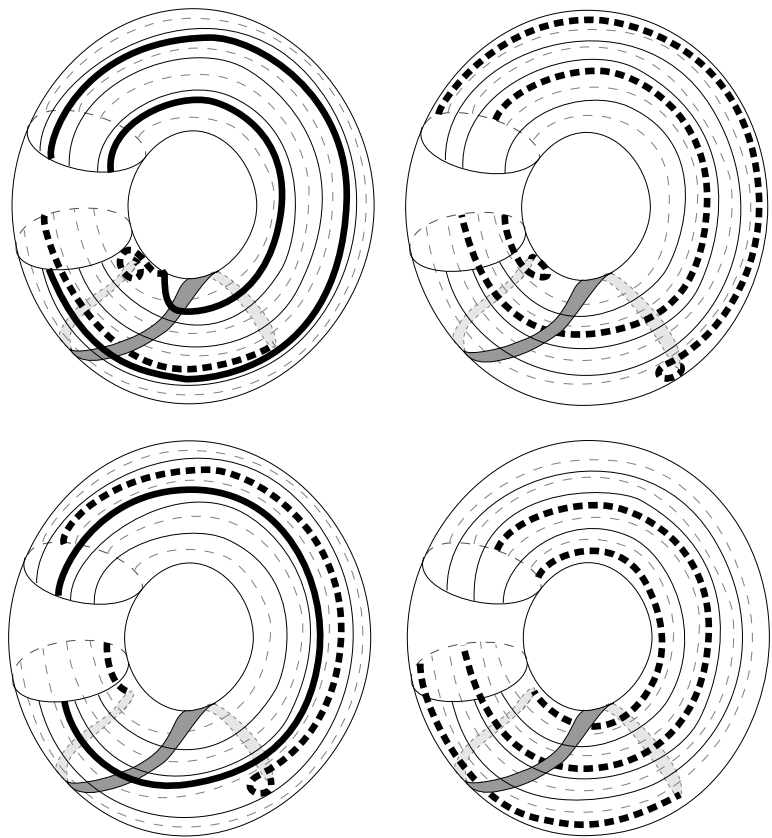}
\caption{The ordinary fiber}
\label{f:fiber}
\end{center}
\end{figure}

Now we consider the portion of the ordinary fiber than lies in the boundary of the handle which contains the tangle $T(p,q)_r$ (from the definition of the twisted torus knot).  The thick black curve in Figure \ref{f:fiber} illustrates the portion of the ordinary fiber curve that lies in this part of $\partial H$, for each of the four cases.  Note that the long, thin gray disk in each case represents the attaching disk $D_r$, after performing the isotopy from Figure \ref{f:long_baguette}.  The blank annulus in each picture represents a $(p,q)$ torus braid.  The light curve is the twisted torus knot $K(p,q,r,1,\epsilon)$.

\begin{figure}[ht!]
\begin{center}
\includegraphics[bb=26 25 203 201]{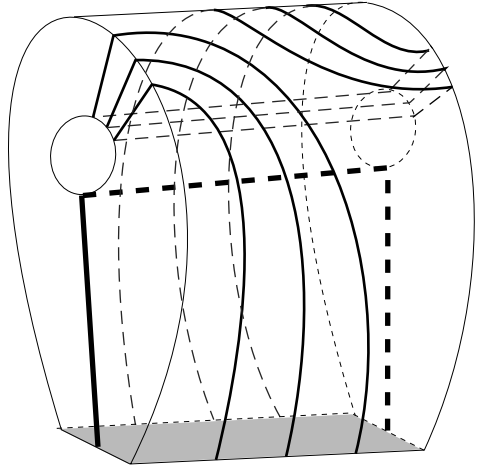}
\caption{Crossing the ``snake"}
\label{f:snake}
\end{center}
\end{figure}

The four representative cases in Figure \ref{f:fiber} are arranged so as to correspond with Table \ref{tb:homology}.    In the top two cases in Figure \ref{f:fiber},  $r=10$, while $r=11$ in the bottom two.  Each picture has $p=7$.  The left two pictures have $\tilde q=2$ while the right two have $\tilde q = -2$.  The top two pictures satisfy $\bar r = p - 2 \hat q$, while the bottom two have $\bar r = 2 \hat q$.

Now we must describe the behavior of the ordinary fiber curve (the thick dark curve) when it goes over the other handle (which is depicted in Figure \ref{f:snake} and contains the tangle $rT(1,\epsilon)$).  There are two ways that the ordinary fiber crosses the handle in Figure \ref{f:snake}; we describe them below.  Since the attaching disk $D_r$ for the handle is the gray, snake-like disk in Figure \ref{f:fiber}, we refer to it as the snake. Note that the attaching disk $D_r$ appears in Figure \ref{f:snake} as the gray rectangular base of the handle.
In each of the four pictures in Figure \ref{f:fiber}, most ordinary fiber seems to cross the snake via a short arc.  In this case, the ordinary fiber crosses the handle as does each of the three parallel arcs in the Figure \ref{f:snake}.  The direction of the twisting depends on $\epsilon$, as usual.  When the ordinary fiber appears in Figure \ref{f:fiber} to cross the ``snake" the long way (entering at the head, and exiting at the tail), it crosses the handle in an untwisted manner as does the thick black arc that in Figure \ref{f:snake}.  Note that the ordinary fiber is disjoint from $K$ in every case.  This completely describes the ordinary fiber.

To verify that each curve is indeed an ordinary fiber, one may check that each curve represents the desired element (already described) in the fundamental group of $W$.  Since these Seifert manifolds (base $D^2$, 2 critical fibers) have incompressible boundary, the curve must be an ordinary fiber in $W$.  One can then easily calculate the homology class of the ordinary fiber curve in the complementary handlebody $H'$; the results are shown in Table \ref{tb:homology}.

\begin{table}
\begin{center}
\begin{tabular} {|l||l|l|} \hline
\strut & $\tilde q > 0$ & $\tilde q < 0$ \\ \hline \hline
\strut $\bar r = p-\alpha \hat q$ & $(\epsilon \alpha \beta - 1, \alpha \gamma + \beta + 1)$ & 
$(\epsilon \alpha \beta + 1, \alpha \gamma - \beta - 1)$ \\ \hline
\strut $\bar r = \alpha \hat q$ & $(\epsilon \alpha (\beta + 1) + 1,  \alpha \gamma - \beta)$ & 
$(\epsilon \alpha (\beta + 1) - 1, \alpha \gamma + \beta)$\\ \hline
\end{tabular}
\caption{The value of $[f]$ in $H_1(H')$}
\label{tb:homology}
\end{center}
\end{table}

\subsection{Determining the Multiplicity of the Third Critical Fiber}

Having located a curve in $\partial H - N(K)$ that becomes an ordinary fiber in $$W = H \cup_K \text{2-handle}$$
we now calculate the multiplicity of the third critical fiber in the small Seifert-fibered space arising from Dehn surgery on the knots in Theorem \ref{t:class}.

Since $K$ is Seifert-fibered with respect to $H$, $W = (H \cup_K \text{2-handle})$ is a Seifert-fibered space over the disk with two critical fibers.  Since $K$ is also primitive with respect to $H'$, we know that $W' \approx D^2 \times S^1$ is the other piece of the manifold obtained by Dehn surgery at the surface slope $\gamma$ (recall Lemma \ref{l:surfslope}).  Let $f$ be a curve in $\partial W$ that is an ordinary fiber of $W$.  If non-zero, the intersection number in $W'$ of the inclusion of $f$ and a meridian disk for $W'$ is the multiplicity $\mu_3$ of the third critical fiber of the Dehn surgery manifold $M(K,\gamma)$.  If this intersection number equals zero then the manifold is a connected sum of lens spaces, as observed earlier.  Our goal, then, is to calculate this intersection number.

Now let $f$ be a curve in $\partial H - N(K)$ that becomes an ordinary fiber under inclusion in $W$.  Since $K$ is primitive, we may choose a disk system $(\Delta _1, \Delta _2)$ for $H'$  and a curve $c$ in $\partial H'$ such that $|\partial \Delta _1 \cap K| = 1$, $|\partial \Delta _2 \cap c| = 1$, $\partial \Delta _2 \cap K = \emptyset$, and  $\partial \Delta _1 \cap c = \emptyset$. Then $\Delta _2$ will become a meridian disk in $W'$, and the algebraic intersection number in $H'$ of $f$ and $\Delta _2$ will equal their algebraic intersection number after inclusion in $W'$.  

Thus to calculate the multiplicity of the third critical fiber, we need only calculate the algebraic intersection number in $H'$ of $f$ and $\Delta _2$. Write $[f]=(f_1, f_2)$, $[K] = (k_1, k_2)$, and $[c] = (c_1, c_2)$ with respect to the usual basis for $H_1(H')$ (i.e.\ the basis arising from $D_1'$ and $D_2'$ in Section \ref{ss:calculatingw}). $[K]$ and $[c]$ also form a basis for $H_1(H')$, so we may write $[f] = s[K] + t[c]$ for some integers $s$ and $t$.  Clearly the algebraic intersection number of $f$ and $\Delta _2$ equals $t$, so we must calculate the second coordinate of $[f]$ when expressed in the $[K], [c]$ basis.  Since $[K]$ and $[c]$ form a basis for $\mathbb{Z}\oplus \mathbb{Z}$, we may assume that $k_1c_2 - k_2c_1 = 1$.  The matrix which changes basis from the usual one to the $[K], [c]$ basis is 
$$\begin{bmatrix} c_2 & -c_1\\ -k_2 & k_1 \end{bmatrix}.$$
Now we apply this change of basis to $(f_1, f_2)$ and observe that the second coordinate is $-k_2f_1 + k_1f_2$.  In summary, we have the following convenient expression for the third multiplicity:
$$ \mu_3 = \begin{vmatrix} k_1 & k_2\\ f_1 & f_2 \end{vmatrix}.$$

In the last subsection we calculated $[f] = (f_1, f_2)$ for any middle Seifert-fibered twisted torus knot. It is easy to see that $[K(p,q,r,1, \epsilon)] = (\epsilon r, q)$ with respect to the usual basis for $H_1(H')$.  Thus, we have all the necessary ingredients to calculate $\mu_3$ for each of the primitive/Seifert-fibered knots described in Theorem \ref{t:class}.  The results of this calculation are given below.

\begin{thm}
\label{t:muls}
The following are the multiplicities $(\mu_1, \mu_2, \mu_3)$of the small Seifert-fibered manifolds arising from Dehn surgery at the surface framing for each of the knots in Theorem \ref{t:class} of this section.
\begin{center}
\begin{tabular} {|l||l|} \hline
\strut  & $(\mu_1, \mu_2, \mu_3)$ \\ \hline \hline
\strut1 & $(2, 2q-p, p+(\epsilon - 2)q)$ \\ \hline
\strut2 & $(k, p-kq, p-(k-\epsilon )q)$ \\ \hline
\strut3 & $(m, l+\delta, m(l-\epsilon \delta)+\delta)$ \\ \hline
\strut4 & $(m, l+1, (-l+n(lm+l+1))(1+\epsilon n - \epsilon) + \epsilon m(nl - l -n))$ \\ \hline
\strut5 & $(m-1, l-1, (m+n(lm-1))(\epsilon n+1+\epsilon) - \epsilon (ln+1) )$ \\ \hline
\end{tabular}
\end{center}
\end{thm}

\section{The ubiquity of SSFS surgeries}

\label{s:ubiquity}
Now we focus on the knots of the second type in Theorem \ref{t:class} (i.e.\ $K(p,q,p-kq,1,\pm 1)$).  Using these knots, we will show that, for any triple of integers $(\mu_1, \mu_2, \mu_3)$ with $(\mu_1,\mu_2)=1$, there is a non-torus knot with a SSFS surgery realizing those multiplicities.  As noted earlier, many of these knots are known to be hyperbolic, and we expect that all knots discussed in this section are hyperbolic.  Note that, for homological reasons, only SSF spaces with $\gcd (\mu_1, \mu_2, \mu_3)$ $=\, 1$ can be realized by surgery on a knot.

The multiplicities of the critical fibers of a small Seifert-fibered space do not determine its homeomorphism class.  However, the determination of which triples of multiplicities are realized does provide a rough measure of which SSF manifolds arise via Dehn surgery on knots.

\subsection{The knots $K(p,q,p-kq,1,1)$}

Here we study those knots in Theorem \ref{t:class} with positive twisting.  While the results in this section are not used to prove the ubiquity theorem mentioned above, they illuminate some features of an interesting class of twisted torus knots.  First we show that the knots $K(p,q,p-kq,1,1)$ as in Theorem \ref{t:class} are not torus knots.  To prove this, we use Lemma \ref{l:nontorus}, whose proof is somewhat surprising.  Lastly, we specify exactly which triples of multiplicities arise by Dehn surgery at the surface slope (namely $pq + (p-kq)^2$) for these knots.  As mentioned above, we conjecture that all of the knots $K(p,q,p-kq,1,1)$ in Theorem \ref{t:class} are, in fact, hyperbolic.

First we show that the knots are non-torus.  If $m/1$ surgery on an $(a,b)$ torus yields a small Seifert manifold, then the multiplicities are $(a, b, |ab-m|)$ \cite{moser:1971}.  Thus, if a knot has an \emph{integral} small Seifert-fibered Dehn surgery for which the sum or difference of the slope and the product of any two coprime multiplicities does not equal the remaining multiplicity (up to sign), then the knot is non-torus.

Both torus knots and the knots $K(p,q,p-kq,1,1)$ are closed positive braids; thus they are fibered, and their fiber genus is easily calculated (\cite{stallings:1978}).  In particular, for a knot which is presented as a closed positive braid, $-\chi (F) = c-s$, where $F$ is the fiber surface, $c$ is the number of crossings, and $s$ is the number of strands of the braid.  We will use this together with the first idea to show that the aforementioned knots are non-torus.

\begin{lem}
\label{l:nontorus}
Let $\kappa$ be a fibered knot with fiber surface $F$.  Suppose that $\kappa$ has an integral SSF surgery with slope $m/1$ and multiplicities $(\mu_1,\mu_2,\mu_3)$.  Then
$$|m| - \mu_1 - \mu_2 - \mu_3 + \chi (F) > 0$$
implies that $\kappa$ is not a torus knot.
\end{lem}

\proof
Note that for each $i$, $\mu_i > 0$.  For convenience, assume $m \ge 0$.  If $m<0$ then we apply the following argument to $-m$ surgery on the mirror image of $K$.

If $\kappa$ were a torus knot, then by the observation above it must be a $(\mu_{i_1}, \pm \mu_{i_2})$ torus knot satisfying
\begin{equation}  \label{e:nont1}
\mu_{i_3} = \pm m \pm \mu_{i_1} \mu_{i_2}
\end{equation}   
for some cyclic permutation $(i_1,i_2,i_3)$ of $(1,2,3)$.  Also, a $(\mu_{i_1}, \pm \mu_{i_2})$ torus knot has a fiber surface $F$ satisfying 
\begin{equation}  \label{e:nont2}
-\chi (F) = \mu_{i_1}\mu_{i_2} - \mu_{i_1} - \mu_{i_2}.
\end{equation}   

To show that $\kappa$ is not a torus knot it is enough to show that equations~\ref{e:nont1} and \ref{e:nont2} have no simultaneous solution for any choice of $(i_1,i_2,i_3)$.  On the surface, it appears that the choice of signs and cyclic permutations might require consideration of twelve such systems of equations.

Fix $(i_1,i_2,i_3)$.  Because $\mu_1$, $\mu_2$, $\mu_3$, and $m$ are all positive, $m - \mu_{i_1}\mu_{i_2} > \mu_{i_3}$ implies that equation \ref{e:nont1} above has no solution for \emph{any} choice of $\pm$ signs.  Thus, it is enough to show that $$m - \mu_{i_1}\mu_{i_2} - \mu_{i_3} > 0$$ for each choice of $(i_1,i_2,i_3)$.

We rewrite equation \ref{e:nont2} as 
$$ \mu_{i_1}\mu_{i_2} - \mu_{i_1} - \mu_{i_2} + \chi (F) =0$$
and add it to the desired inequality, obtaining
$$m - \mu_{i_1} - \mu_{i_2}- \mu_{i_3} + \chi (F) > 0.$$

Notice that, surprisingly, this inequality is symmetric in $i_1$, $i_2$, and $i_3$.  This symmetry implies that one inequality suffices to show that $\kappa$ is not a torus knot.
\endproof

\begin{thm}
\label{t:nontorus}
The twisted torus knots $K(p,q,p-kq,1,1)$ with $1 \le k < p/q$ and $1 < q < p/2$ are not torus knots.
\end{thm}

\proof
Using the positive braid picture of the knots $K(p,q,p-kq,1,1)$, we find that 
$$-\chi(F) = pq-p-q+(p-kq)(p-kq-1),$$
where $F$ is the fiber surface for $K$.

By earlier calculations, we have that $(\mu_1, \mu_2, \mu_3) = (k,p-kq,p-kq+q)$
and $m = pq + (p-kq)^2$.  Note that $m \ge 0$.  Hence, 
$$m - \mu_1 - \mu_2 - \mu_3 + \chi (F) = k(q-1).$$
Since $k \ge 1$ and $q>1$ this quantity is positive, so the lemma implies that $K$ can not be a torus knot.
\endproof

\begin{rems}
Lemma \ref{l:nontorus} can be more powerful than using equation \ref{e:nont1} or equation \ref{e:nont2} separately.  In particular, the slope criterion and the genus criterion each fail individually for some of the knots in Theorem \ref{t:nontorus}.

For example, $K(25,2,5,1,1)$ has a slope $75/1$ SSF Dehn surgery with multiplicities $(10,5,7)$.  The $(10,7)$ torus knot also has a slope $75/1$ surgery with these multiplicities.

An example of the second type is $K(33,2,5,1,1)$, which has a SSF Dehn surgery with multiplicities $(14,5,7)$ and has fiber genus $51$.  The $(14,5)$ torus knot also has fiber genus $51$ and a $(14,5,7)$ SSF Dehn surgery.
\end{rems}

Now we describe exactly which triples of multiplicities arise from $pq+~(p-kq)^2$ surgery on $K(p,q,p-kq,1,1)$.

\begin{thm}  
\label{t:preplethora}
The triples of multiplicities of the SSF manifolds obtained by $pq + (p-kq)^2$ surgery on $K(p,q,p-kq,1,1)$ are exactly those $(\mu_1, \mu_2, \mu_3)$ with $(\mu_1, \mu_2) =1$ and  $|\mu_1-\mu_2| > 1$.
\end{thm}

\proof
First we show that all such triples arise.  We may assume that $\mu_1 > \mu_2$.  The following equations describe knots $K(p,q,p-kq,1,1)$ which realize the desired multiplicities:
\begin{align*} 
q    &= \mu_1-\mu_2\\
p   &= \mu_3(\mu_1 - \mu_2) + \mu_2\\
k    &=\mu_3 
\end{align*}
Since it is easy to verify that $(p,q)=1$ and $q$ and $k$ are in the required range, we conclude that $K(p,q,p-kq,1,1)$ is indeed a knot of the second type in Theorem \ref{t:class}. Using the table in Theorem \ref{t:muls}, one can verify directly that the multiplicities are as desired.

Now we show that no other triples arise.  By Theorem \ref{t:muls}, two of the multiplicities are $p-kq+q$ and $p-kq$.  These integers are coprime (since $p$ and $q$ are coprime) and differ by $q$, which is greater than 1.
\endproof

\subsection{The knots $K(p,q,p-kq,1,-1)$}

In this section we focus on the knots $K(p,q,p-kq,1,-1)$, i.e.\ the second type in Theorem \ref{t:class} with negative twisting.  We use these to prove the following theorem.  As mentioned earlier, we expect that these knots are all hyperbolic, and many are known to be hyperbolic.

\begin{thm}  
\label{t:plethora}
For any triple of multiplicities $(\mu_1, \mu_2, \mu_3)$ with $(\mu_1, \mu_2) =1$, there is a non-torus knot with a small Seifert-fibered Dehn surgery with the prescribed multiplicities.
\end{thm}

\proof
The first step in the proof is almost the same as for Theorem \ref{t:preplethora}.  Let $(\mu_1,\mu_2,\mu_3)$ be a triple of (positive) integers as above.  Set 
\begin{align*} 
q    &= \mu_1+\mu_2\\
p    &= \mu_3(\mu_1 + \mu_2) + \mu_2\\
k    &= \mu_3 
\end{align*}
Then it is easy to verify that the knot $K(p,q,p-kq,1,-1)$ is indeed a knot of the second type in Theorem \ref{t:class}.  By Theorem \ref{t:muls}, $pq - (p-kq)^2$ surgery yields a $(\mu_1,\mu_2,\mu_3)$ SSFS.

Now we show that each of these knots is non-torus.  The surgery slope $m$ is equal to $\mu_3(\mu_1+\mu_2)^2 + \mu_1\mu_2$.  It is very easy to check that $m - \mu_{i_1}\mu_{i_2} > \mu_{i_3}$ for each cyclic permutation $(i_1,i_2,i_3)$ of $(1,2,3)$.  It follows from the slope criterion for torus knot surgeries mentioned earlier that such a SSFS cannot arise from \emph{integral} surgery on a torus knot.

This implies that the knots constructed above (a subclass of the knots $K(p,q,$
$p-kq,1,-1)$) are not torus knots.
\endproof

The finite $\pi_1$ SSF manifolds comprise all known finite-$\pi_1$ $3$-manifolds (allowing trivial fibers); these manifolds are the well-known spherical space forms.  The corresponding spherical triples satisfy $\Sigma_{i=1}^{3} \frac{1}{\mu_i} > 1$.  Thus, the spherical triples (with no trivial fibers) are $(2,3,3)$, $(2,3,4)$, $(2,3,5)$, and $(2,2,n)$.  For homological reasons, $n$ must be odd for a $(2,2,n)$ SSFS to arise from surgery on a knot.  The following Corollary follows immediately from Theorem \ref{t:plethora}.

\begin{cor}
There are non-torus knots in $S^3$ with SSF Dehn surgeries realizing all possible spherical triples of multiplicities, \emph{i.e.} $(2,3,3)$, $(2,3,4)$, $(2,3,5)$, and $(2,2,n)$ for odd $n$.
\end{cor}

As mentioned in the introduction, Theorem \ref{t:plethora} has also been obtained by Miyazaki and Motegi using the knots in \cite{mm:1999}.  Moreover, they show that the knots they construct are indeed hyperbolic.  Many of the spherical triples can be realized using the examples mentioned earlier of Eudave-Mu\~noz, Bleiler and Hodgson, and Boyer and Zhang.

\begin{note}
Are all triples of multiplicities $(\mu_1, \mu_2, \mu_3)$ with $\gcd(\mu_1, \mu_2, \mu_3) = 1$ realized by Dehn surgery on a hyperbolic knot? Looking at Theorem \ref{t:plethora}, the only ones that might be missing are triples with each $\gcd (\mu_i,\mu_j) \ne 1$, but $\gcd(\mu_1,\mu_2,\mu_3) = 1$.  Such a small Seifert-fibered space could not be realized by any surgery on a torus knot.  Examples of this type do arise from knots of type $4$ and $5$ in Theorem \ref{t:class}.  For example, choosing $l=5$, $m=14$, $n=2$, and $\epsilon=-1$ gives a knot with a $(14,6,105)$ SSFS surgery.  Such triples also arise from the Eudave-Mu\~noz examples.
\end{note}

\Addresses\recd

\end{document}